\documentclass[twocolumn]{autart}    

\usepackage{amsmath,amssymb,amsfonts,bm}
\usepackage{graphicx}
\usepackage{algorithm}
\usepackage{algorithmicx,algpseudocode}
\usepackage{cite}
\usepackage{color}
\let\theoremstyle\relax

\usepackage{amsthm}
\usepackage{cuted}
\usepackage[mathscr]{eucal}
\theoremstyle{plain}
\newtheorem{theorem}{Theorem}[section]
\newtheorem{lemma}[theorem]{Lemma}

\newtheorem{proposition}[theorem]{Proposition}
\theoremstyle{definition}

\newtheorem{problem}[theorem]{Problem}
\theoremstyle{remark}

\usepackage{wrapfig}
\newcommand{\x}{\mathrm{x}}
\renewcommand{\u}{\mathrm{u}}
\usepackage{wrapfig}

\begin{document}

\begin{frontmatter}

\title{General Distribution Steering: A Sub-Optimal Solution by Optimization} 

\author[Guangyu]{Guangyu Wu}\ead{guangyu.wu@ntu.edu.sg},    
\author[{Anders1},{Anders2}]{Anders Lindquist}\ead{alq@kth.se}              

\address[Guangyu]{College of Computing and Data Science, Nanyang Technological University, Singapore}  
\address[Anders1]{School of Artificial Intelligence, Anhui University, Hefei, China}
\address[Anders2]{Smart Sensor Fusion Research Laboratory, Shanghai Jiao Tong University, Shanghai, China}
          

\begin{abstract}                          
General distribution steering is intrinsically an infinite-dimensional problem, when the continuous distributions to steer are arbitrary. 
We put forward a moment representation of the primal system for control in \cite{wu2023density}. However, the system trajectory was a predetermined one without optimization towards a design criterion, which doesn't always ensure a most satisfactory solution. In this paper, we propose an optimization approach to the general distribution steering problem of the first-order discrete-time linear system, i.e., an optimal control law for the corresponding moment system. The domain of all feasible control inputs is non-convex and has a complex topology. We obtain a subset of it by minimizing a weighted sum of squared integral distances alongside the system trajectory. The feasible domain is then proved convex, and the optimal control problem can be treated as a convex optimization or by exhaustive search, based on the type of the cost function. Algorithms of steering for continuous and discrete distributions are then put forward respectively, by adopting a realization scheme of control inputs. We also provide an explicit advantage of our proposed algorithm by truncated power moments to the prevailing Gaussian Mixture Models. Experiments on different types of cost functions are given to validate the performance of our proposed algorithm. Since the moment system is a dimension-reduced counterpart of the primal system, we call this solution a sub-optimal one to the primal general distribution steering problem.
\end{abstract}

\end{frontmatter}

\section{Introduction}
In this paper, we consider the problem of steering the distribution of the state where the system dynamics are governed by a discrete-time stable first-order linear stochastic difference equation. The linear dynamics of the system read
\begin{equation}
x(k+1)=a(k) x(k)+u(k) .
\label{uniequation}
\end{equation}
Here stability refers to the fact that the equilibrium $x(k)=0$ is asymptotically and exponentially stable. Since the system is stable and we assume $a(k)$ to be positive, we have $a(k) \in(0,1)$. The control input to the system at time step $k$ is defined as $u(k)$, and $x(k)$ is its state. 

Distribution steering is widely prevalent across various control domains. In some scenarios, it plays a vital role in characterizing the uncertainty of system states in certain scenarios \cite{chen2015optimal, chen2015optimal2, chen2016optimal, balci2020covariance, okamoto2018optimal}. Brockett also emphasized the significance of considering the distribution of the system state in optimal control \cite{brockett2007optimal}. On the other hand, distribution steering has become a core research area of swarm robotics \cite{dorigo2020reflections, dorigo2021swarm} recently. Mean-field model is widely adopted for controlling a vast group of swarm robots \cite{liu2018mean, deshmukh2018mean, elamvazhuthi2019mean}. As the number of agents approaches infinity, it becomes feasible to approximate the linear forward equation of each agent in a whole with a forward equation with parameters being continuous probability distributions. The resultant equation, referred to as the mean-field model, is formulated based on a set of probability densities that govern the probability of an agent occupying a particular state at a given time. In cases where the swarm consists of a large number of agents, this approximation remains accurate provided that all agents adhere to identical control laws (i.e., each agent is homogeneous), and each agent's control law is solely influenced by its own state or the local density of the swarm, rather than being contingent upon the identities of other agents. Distribution steering of discrete distributions in the form of occupation measures, which serves as a compact representation of the agents, is also considered in \cite{wu2022group}.

A definition of discrete-time general distribution steering treated in this paper is given as follows. Provided with an initial random variable $x(0)$, the distribution steering problem amounts to choosing a sequence of random variables $(u(0), u(1), \cdots, u(K-1))$, so that the probability distribution $\chi_0$ of $x(0)$ is transferred to the distribution $\chi_K$ of $x(K)$ at some future time $K$. For the general distribution steering problem, as is proposed in \cite{sivaramakrishnan2022distribution}, the distributions of all $x(k)$ and $u(k)$ are all arbitrary, which are not assumed to fall within specific classes of function (e.g. Gaussian). The distributions can be either continuous or discrete. We note that when the distributions are continuous, the general distribution steering problem is intrinsically an infinite-dimensional one (both the functional space and the parameter space of the distributions are infinite-dimensional).

The distribution steering problem has a history of decades \cite{hotz1987covariance, collins1985covariance, collins1987theory, hsieh1990all, xu1992improved}, and has recently re-emerged as a hot topic in control theory and engineering, owing to its theoretical elegance and practical applications in areas such as swarm robotics and flow modeling. Roughly speaking, there are two main lines of research on the distribution steering problem. The first line considers settings where the system has no intrinsic dynamics in the absence of control inputs, namely the control input alone drives the evolution of the system state over time. This formulation is particularly common in swarm robotics, where the system state represents the positions of robots, and the control inputs determine their collective motion. Zheng, Han, and Lin \cite{zheng2021transporting, zheng2020pde, zheng2021distributed} used mean-field partial differential equations, specifically the Fokker–Planck equation, to model the swarm and control the mean-field density of the velocity field. Biswal, Elamvazhuthi, and Berman \cite{elamvazhuthi2017mean, biswal2021decentralized} approached the problem by stabilizing the Kolmogorov forward equation, which describes the mean-field behavior of the swarm. Caluya and Halder \cite{caluya2021wasserstein} proposed Wasserstein proximal algorithms for solving this problem in the framework of the Schr\"odinger bridge.

For the other line of research, people consider the system dynamics of the group of agents to be steered. This type of distribution steering is more general than the first type, however is inevitably more difficult. As a tradeoff, the distribution of the agents are assumed to fall within specific classes of functions to ensure the solvability of the problem. A most widely considered distribution is the Gaussian, which ensures a closed form of solution to the distribution steering problem. The distribution steering problem is then reduced to steering the statistics of the distribution. For the Gaussians, the task is to steer the mean and variance of the distributions, which is called ``covariance steering'' in the literature. Pioneering results for covariance steering include \cite{okamoto2018optimal, okamoto2019optimal} by Okamoto and Tsiotras, \cite{liu2023optimal, liu2022optimal1, liu2022optimal2} by Liu and Tsiotras, \cite{yin2022trajectory} by Yin, Zhang, Theodorou and Tsiotras and \cite{saravanos2022distributed} by Saravanos, Balci, Bakolas and Theodorou. Moreover, Sivaramakrishnan, Pilipovsky, Oishi and Tsiotras \cite{sivaramakrishnan2022distribution} proposed to treat the non-Gaussian distribution steering problem by characteristic functions, which was one of the earliest attempts for the general distribution steering. For the continuous-time linear systems, Chen, Georgiou and Pavon have proposed fundamental results using the Schr\"odinger Bridge strategy for Gaussian distributions \cite{chen2015optimal, chen2015optimal2} and other types of distributions \cite{chen2018optimal}. Caluya and Halder \cite{caluya2021reflected} extended the results to nonlinear continuous-time systems and hard state constraints. Moreover, Sinigaglia, Manzoni, Braghin and Berman \cite{sinigaglia2022robust} put forward a robust optimal density control of robotic swarms. 

The findings presented above, along with numerous others, have significantly contributed to the distribution steering problem. However, conventional assumptions have typically constrained distributions to specific function classes. In many practical scenarios, such as controlling a collective of agents, as will be addressed in subsequent sections of this paper, assuming the agents' distribution to be Gaussian is always improper. On the other hand, to the best of our knowledge, a complete solution for the distribution steering problem considering discrete-time linear systems, where both initial and terminal distributions are non-Gaussian, hasn't been proposed. Furthermore, given the infinite-dimensional nature of the general distribution steering problem, the presence of error in the solution is inevitable, rendering the problem open and consequently non-trivial. How to quantify the error of steering should also be taken into consideration.

Let's turn our focus to another way of characterizing the probability distribution. In the probability theory, we know that a distribution function can be uniquely determined by its full power moment sequence \cite{Shiryaev2016}. The previously proposed distribution steering problem aims to regulate the system state represented as a probability distribution. When the distribution is assumed to be Lebesgue integrable only, the problem becomes uncountably infinite-dimensional, which is generally intractable. Steering the system by controlling the full power moment sequence, rather than directly manipulating the distribution of the system state, reduces the problem to a countably infinite-dimensional one, which is still not practically feasible. However, by properly truncating the first several terms of the power moment sequence for characterizing the density of the system state \cite{byrnes2006generalized, wu2023non}, the problem is now steering a truncated power moment sequence to another, which is finite-dimensional and tractable. It is not the first time in the literature that the power moments are used for control purposes. Jasour, together with Lagoa \cite{jasour2014reconstruction} proposed to reconstruct the support of a measure from its moments. It works well for the uniform distributions. Partly based on this result, He, Wang and Williams \cite{jasour2021moment} addressed the problem of uncertainty propagation through the control of nonlinear stochastic dynamical systems. In our previous result \cite{wu2022group}, we proposed to give a reduced-order counterpart of the primal system by the power moments, and to perform controls on the moment system. However, the control law in that paper was predetermined, similar to the treatment in \cite{sivaramakrishnan2022distribution}. We were not able to design the control inputs by desired criteria through optimization in the manner of the conventional optimal control. 

This paper delves into the general distribution steering of first-order discrete-time linear stochastic systems, where the specified initial and terminal distributions are arbitrary, requiring only the existence of the first several orders of finite power moments. We approach this investigation through optimization. The paper's structure is outlined as follows. In Section 2, we introduce a moment-based representation of the primal discrete-time linear system and formulate the distribution steering problem using the moment system. We also explore the controllability of the moment system. Section 3 presents a optimization framework for controlling the moment system. Due to the necessity for positive definiteness in the Hankel matrices of both the moments of control inputs and system states, the domain of feasible moments for control inputs, given the desired terminal moments of the system state, is non-convex and topologically complex. We propose a domain for optimization and establish its convexity. We prove that this choice of the state of the moment system essentially minimizes a weighted sum of squared integral distances. Choices of cost functions are provided in Section 4, and corresponding optimization approaches are proposed. In Section 5, we employ a distribution parametrization algorithm proposed in our previous work \cite{wu2022non1} to realize control inputs as analytic functions using power moments obtained from the proposed control scheme, which is a convex optimization algorithm with the existence and uniqueness of solution proved. Section 6 introduces algorithms for two types of distribution steering problems: continuous distribution steering and discrete distribution steering. In Section 7, we compare our proposed algorithm, based on truncated power moments, with Gaussian Mixture Models—an existing dimensionality reduction technique for general distribution steering—and provide a detailed explanation of the advantages of our approach. For simulation purposes, we consider two typical distribution steering problems in Section 8. The first problem involves separating a group of agents into smaller groups, and the second aims to steer agents from two groups to desired terminal groups. Numerical examples demonstrate the performance of our proposed algorithms under different types of cost functions. It provides a new perspective to distributional control as Brockett advocated in \cite{brockett2007optimal}.

\section{A moment formulation of the primal problem}
In this section we treat the distribution steering problem formulated in Section 1. In the conventional distribution steering schemes, feedback control is usually adopted, i.e., the control inputs are usually chosen as functions of system states. It is easy to verify that by the feedback control, namely $u(k)$ being a function of $x(k)$, the terminal distribution $\chi_{K}$ is the initial distribution $\chi_{0}$ after alternations of location and scale. Therefore, feedback control law is not suitable for the general distribution steering problem we consider, where the initial and terminal distributions are generally from different function classes. 

Distinct from the conventional control strategies, we extend the control input to a random variable rather than a function of the system state. Since the feedback control $u(k)$ is a function of the random variable $x(k)$, it is indeed a random variable. By extending $u(k)$ to a general random variable, which is not a function of $x(k)$, we will be able to steer a distribution to one having a different function type, in the following parts of this paper. However by doing this, it is still not always possible to obtain a closed-form solution to this problem. If the distributions are not assumed to fall within certain specific classes, the problem is intrinsically infinite-dimensional. Define the distribution of the control $u(k)$ as $\nu_{k}(\u)$. We further assume the system states $x(k)$ and the control inputs $u(k)$ being independent. This assumption is not the first time in the literature, which has already been used in \cite{van2021control} for treatments of stochastic control systems. By this assumption the distribution of $x(k+1)$ can be written as

\begin{equation}
\begin{aligned}
    \chi_{k+1}(s) & = \int_{\mathbb{R}} \chi_{k}\left(\frac{\varepsilon}{a(k)}\right) \nu_{k}\left(s-\varepsilon\right) d\varepsilon\\
    & = \left(\chi_k\left(\frac{s}{a(k)}\right) * \nu_k(s)\right)(s).
\end{aligned}
\label{qk1}
\end{equation}

For the distribution steering problem, an analytic solution of $\chi_{k+1}(s)$ in \eqref{qk1} is necessary. However, except for limited classes of functions such as Gaussian distributions and trigonometric functions, this isn't possible in general. When $x(k)$ and $u(k)$ are not independent, the problem is even more complicated. This is the main reason that in previous results which have a similar problem setting, the examples have almost always Gaussian or trigonometric densities. This severely limits the use of these results in real applications.

A similar problem exists in non-Gaussian Bayesian filtering. In our previous results \cite{wu2023non}, we proposed a method of using the truncated power moments to reduce the dimension of this problem, mainly for characterizing the macroscopic properties of the distributions. This strategy can also be found in \cite{byrnes2003convex, georgiou2003kullback}, which turns the problem we treat to a tractable truncated moment problem.

Denote the expectation operator as $\mathbb{E}[\cdot]$. By the system equation \eqref{uniequation}, the power moments of the states up to order $2n$ are written as
\begin{equation}
\mathbb{E}\left[ x^{l}(k+1) \right] = \sum_{j=0}^{l}\binom{l}{j}a^{j}(k)\mathbb{E}\left[ x^{j}(k)u^{l-j}(k) \right].
\label{MomentCal}
\end{equation}

We note that it is difficult to treat the term $\mathbb{E}\left[ x^{j}(k)u^{l-j}(k) \right]$. However, we note that if $x(k)$ and $u(k)$ are independent, i.e., $\mathbb{E}\left[ x^{j}(k)u^{l-j}(k) \right] = \mathbb{E}\left[ x^{j}(k) \right]\mathbb{E}\left[ u^{i-j}(k) \right]$, the dynamics of the moments can be written as a linear matrix equation
\begin{equation}
    \mathfrak{X}(k+1) = \mathfrak{A}(\mathfrak{U}(k))\mathfrak{X}(k)+\mathfrak{U}(k)
\label{momentsystem}
\end{equation}
where the state vector is composed of the power moment terms up to order $2n$, i.e.,
\begin{equation}
\mathfrak{X}(k) = \begin{bmatrix}
\mathbb{E}[x(k)] & \mathbb{E}[x^{2}(k)] & \cdots & \mathbb{E}[x^{2n}(k)]
\end{bmatrix}^{T},
\label{XK}
\end{equation}
and the input vector is written as
\begin{equation}
\mathfrak{U}(k) = \begin{bmatrix}
\mathbb{E}[u(k)] & \mathbb{E}[u^{2}(k)] & \cdots & \mathbb{E}[u^{2n}(k)] 
\end{bmatrix}^{T} .
\label{UK}
\end{equation}
Here 
\begin{equation}
\mathbb{E}\left[ x^{l}(k) \right] = \int_{\mathbb{R}}\x^{l}\chi_{k}(\x)d\x
\label{xlK}
\end{equation}
and
$$
\mathbb{E}\left[ x^{j}(k)u^{l-j}(k) \right] = \int_{\mathbb{R}}\x^{j}\chi_{k}(\x)d\x \int_{\mathbb{R}}\u^{l-j}\nu_{k}(\u)d\u.
$$
for $l \in \mathbb{N}_{0}$ ($\mathbb{N}_{0}$ denotes the set of all nonnegative integers), $l \leqslant 2n$.
Similarly we have
\begin{equation}
\mathbb{E}\left[ u^{l}(k) \right] = \int_{\mathbb{R}}\u^{l}\nu_{k}(\u)d\u .
\label{ulK}
\end{equation}
The matrix $\mathfrak{A}(\mathfrak{U}(k))$ in the system \eqref{momentsystem} can then be written as \eqref{longeq0}.

\begin{figure*}[t]
\begin{equation}
\mathfrak{A}(\mathfrak{U}(k))
= \begin{bmatrix}
a(k) & 0 & 0 & \cdots & 0\\ 
2a(k)\mathbb{E}[u(k)] & a^{2}(k) & 0 & \cdots & 0\\ 
3a(k)\mathbb{E}[u^{2}(k)] & 3a^{2}(k)\mathbb{E}[u(k)] & a^{3}(k) & \cdots & 0\\ 
\vdots & \vdots & \vdots & \ddots\\ 
\binom{2n}{1}a(k)\mathbb{E}[u^{2n-1}(k)] & \binom{2n}{2}a^{2}(k)\mathbb{E}[u^{2n-2}(k)] & \binom{2n}{3}a^{3}(k)\mathbb{E}[u^{2n-3}(k)] &  & a^{2n}(k)
\end{bmatrix}.
\label{longeq0}
\end{equation}
\hrulefill
\vspace*{4pt}
\end{figure*}

By using the truncated power moments to characterize the dynamics of system \eqref{uniequation} where $x(k)$ and $u(k)$ are random variables, we shall reformulate the control problem as steering the power moments of the $x(k)$ and $u(k)$. System \eqref{momentsystem} is called the moment system corresponding to system \eqref{uniequation}. The power moment steering problem is then formulated as follows.

\begin{problem}
The dynamics of the moment system is
$$
    \mathfrak{X}(k+1) = \mathfrak{A}(\mathfrak{U}(k))\mathfrak{X}(k)+\mathfrak{U}(k)
$$
where $\mathfrak{X}(k), \mathfrak{U}(k)$ are obtained by \eqref{xlK} and \eqref{ulK}. Given an arbitrary initial distribution $\chi_{0}(\x)$ (the first $2n$ orders of power moments exist) and terminal power moments $\{\mu_{i}\}_{i = 1:2n}$, determine the control sequence $$\left( u(0), \cdots, u(K-1)\right)$$ so that the first $2n$ order power moments of the terminal distribution are identical to those specified, i.e.,
\begin{equation}
\mathfrak{X}(K) = \int_{\mathbb{R}} \x^{l} \chi_{K}(\x) d\x = \mu_{l}
\label{ExTl}
\end{equation}
for $l = 1, \cdots, 2n$.
\end{problem}

However for the moment system to control, there remains to design control laws which satisfy 
\begin{equation}
\mathbb{E}\left[ x^{j}(k)u^{l-j}(k) \right] = \mathbb{E}\left[ x^{j}(k) \right]\mathbb{E}\left[ u^{i-j}(k) \right].
\label{Expeq}
\end{equation}

To satisfy \eqref{Expeq}, the control vector is required to be independent of the current state vector. In the conventional feedback control law, this is hardly possible since the control inputs are always functions of the state vectors. However, for distribution steering problems, we note that it is possible to satisfy \eqref{Expeq}, since the control inputs of the primal system, as well as the system states, are probability distributions. For a given system state, by drawing an i.i.d. sample from the probability distribution of the control input, we are able to obtain a control input which is independent of the current system state. By doing this, $x(k)$ and $u(k)$ are independent, i.e., \eqref{Expeq} is satisfied.

Moreover, we note that the control inputs in the moment system are essentially the power moments of the controls to the primal system \eqref{uniequation}. For the univariate random variables, the sufficient and necessary condition of existence is the positive definiteness of the Hankel matrix. The Hankel matrix of $\mathfrak{X}(k)$ reads
$$
\left[ \mathfrak{X}(k) \right]_{H} = \begin{bmatrix}
1 & \mathbb{E}\left[x(k)\right ] & \hdots & \mathbb{E}\left[x^{n}(k)\right ]\\ 
\mathbb{E}\left[x(k)\right ] & \mathbb{E}\left[x^{2}(k)\right ] & \hdots & \mathbb{E}\left[x^{n+1}(k)\right ]\\ 
\vdots & \vdots & \ddots & \\ 
\mathbb{E}\left[x^{n}(k)\right ] & \mathbb{E}\left[x^{n+1}(k)\right ] &  & \mathbb{E}\left[x^{2n}(k)\right ]
\end{bmatrix},
$$
where $\left[ \mathfrak{X}(k) \right]_{H}$ denotes the Hankel matrix. Moreover, we define a subspace of $\mathbb{R}^{2n}$ as $\mathbb{V}^{2n}_{++} := \{ \mathfrak{X} \in \mathbb{R}^{2n} \mid \left[\mathfrak{X}\right]_{H} \succ 0 \}$. Different from the conventional control problems, we confine both $\mathfrak{X}(k)$ and $\mathfrak{U}(k)$ for $k = 0, \cdots, K-1$ to fall within $\mathbb{V}^{2n}_{++}$ to ensure the existence of the corresponding $x(k)$ and $u(k)$. It makes the problem more complicated than usual. Therefore, before we really settle down to treat the control of the moment system \eqref{momentsystem}, we would first like to prove the controllability of it.

\begin{theorem}[Controllability of system \eqref{momentsystem}]
Given system equation \eqref{momentsystem}, there exists a $K$, satisfying $K < +\infty$ and $K \in \mathbb{N}_{0}$, such that an arbitrary initial $\mathfrak{X}(0)$ can be steered to an arbitrary $\mathfrak{X}_{T}$ within $K$ steps.
\label{Theorem31}
\end{theorem}

\begin{proof}
It suffices to prove that there always exists a control sequence $(\mathfrak{U}(0), \cdots, \mathfrak{U}(K - 1))$, which is feasible of steering an arbitrary $\mathfrak{X}(0)$ to an arbitrary $\mathfrak{X}(K)$.

We propose the following control strategy. Before time step $k_{0}$, the system is uncontrolled, i.e., $\mathfrak{U}(k) = 0$ for $k \leqslant k_{0}$. Then we have
$$
\mathfrak{X}(k_{0}) = \mathfrak{A}_{0:k_{0} - 1}\left( 0 \right)\mathfrak{X}(0)
$$
where
$$
\mathfrak{A}_{0:k_{0} - 1}\left( 0 \right) = \begin{bmatrix}
\prod_{k=0}^{k_{0}-1}a(k) &  &  \\
 &  \ddots &  \\
 &  &  \prod_{k=0}^{k_{0}-1}a^{2n}(k) \\
\end{bmatrix}.
$$

We then have
\begin{equation}
\lim_{k_{0} \rightarrow +\infty}\mathfrak{A}_{0:k_{0} - 1}\left( 0 \right)\mathfrak{X}\left( 0 \right) = \mathbf{0}.
\label{k0inf}
\end{equation}

Substituting $k = k_{0}$ into \eqref{momentsystem}, we have
$$
\begin{aligned}
& \mathfrak{X}(k_{0}+1)\\
= & \mathfrak{A}(\mathfrak{U}(k_{0}))\mathfrak{X}(k_{0})+\mathfrak{U}(k_{0})\\
= & \mathfrak{A}(\mathfrak{U}(k_{0}))\mathfrak{A}_{0:k_{0} - 1}\left( 0 \right)\mathfrak{X}(0)+\mathfrak{U}(k_{0}).\\
\end{aligned}
$$

Since $x(k_{0}+1)$ exists, it is obvious that $\left[\mathfrak{X}(k_{0}+1))\right]_{H} \succ 0$. By \eqref{k0inf}, there always exists a $k_{0}<+\infty$ such that
$$
\begin{aligned}
& \mathfrak{U}(k_{0})\\
= & \mathfrak{X}(k_{0}+1) - \mathfrak{A}(\mathfrak{U}(k_{0}))\mathfrak{A}_{0:k_{0} - 1}\left( 0 \right)\mathfrak{X}(0) \in \mathbb{V}^{2n}_{++}.
\end{aligned}
$$

Therefore, there exists at least a control sequence $(0, \cdots, 0, \mathfrak{U}(k_{0}))$, which is feasible of steering an arbitrary initial $\mathfrak{X}(0)$ to an arbitrary $\mathfrak{X}_{T}$ within $K = k_{0}+1$ steps.
\end{proof}

\section{An optimization scheme}

Suppose we are now confronted with the distribution steering problem for system \eqref{momentsystem}, of which the initial moment vector is $\mathfrak{X}(0)$ and the terminal moment vector is $\mathfrak{X}_{T}$ as desired. 

It would be natural for one to consider obtaining the moment vectors of the controls by the following optimization scheme
\begin{equation}
\begin{aligned}
& \operatorname{minimize} f(\mathfrak{U}(0), \cdots, \mathfrak{U}(K-1))\\
\text{s.t.} \quad & \mathfrak{X}(k+1) = \mathfrak{A}(\mathfrak{U}(k))\mathfrak{X}(k)+\mathfrak{U}(k),\\
\quad & \mathfrak{X}(K) = \mathfrak{X}_{T}, \ \mathfrak{U}(k) \in \mathbb{V}^{2n}_{++}.
\end{aligned}
\label{optimization}
\end{equation}
where $f(\cdot)$ is a cost function. Furthermore, by selecting $f(\cdot)$ as a convex function, the optimization problem \eqref{optimization} is convex, given that the following set
$$
\begin{aligned}
\mathcal{U}_{\mathfrak{X}_{T}} := \{ & \left(\mathfrak{U}(0), \cdots, \mathfrak{U}(K-1)\right) \mid \mathfrak{X}(k+1) = \mathfrak{A}(\mathfrak{U}(k))\mathfrak{X}(k)\\
& +\mathfrak{U}(k),
\mathfrak{X}(K) = \mathfrak{X}_{T} \}
\end{aligned}
$$
is convex. However, it is not the case, which will be proved in the following lemma.

\begin{lemma}
The set $\mathcal{U}_{\mathfrak{X}_{T}}$ is not convex, given $K > 1$.
\label{Lemma31}

\end{lemma}
\begin{proof}
Let us assume two sequences 
$$
\left( \mathfrak{U}'(0), \cdots, \mathfrak{U}'(K-1) \right) \in \mathcal{U}_{T}$$
and 
$$
\left( \mathfrak{U}''(0), \cdots, \mathfrak{U}''(K-1) \right) \in \mathcal{U}_{T}.
$$ 
For the set $\mathcal{U}_{T}$ to be convex, we need to have 

$$
\begin{aligned}
& \left( \lambda \mathfrak{U}'(0) + \left( 1 - \lambda \right)\mathfrak{U}''(0), \cdots, \lambda \mathfrak{U}'(K-1) \right.\\
& \left. + \left( 1 - \lambda \right)\mathfrak{U}''(K-1) \right) \in \mathcal{U}_{T}, \quad \forall \lambda \in \left(0, 1\right)
\end{aligned}
$$

Since the two sequences are in the set $\mathcal{U}_{T}$, we have
$$
\mathfrak{X}_{1}(1) = \mathfrak{A}(\mathfrak{U}'(0))\mathfrak{X}(0)+\mathfrak{U}'(0)
$$
and
$$
\mathfrak{X}_{2}(1) = \mathfrak{A}(\mathfrak{U}''(0))\mathfrak{X}(0)+\mathfrak{U}''(0).
$$

By \eqref{momentsystem} we have
\begin{equation}
\begin{aligned}
& \mathfrak{A}\left(\lambda \mathfrak{U}'(0) + \left(1 - \lambda \right) \mathfrak{U}''(0)\right)\mathfrak{X}(0)\\
+ & \lambda\mathfrak{U}'(0)+\left(1 - \lambda\right)\mathfrak{U}''(0)\\
= & \lambda \left(\mathfrak{A} \left(\mathfrak{U}'(0)\right) \mathfrak{X}(0) + \mathfrak{U}'(0)\right)\\
+ & \left(1 - \lambda\right) \left(\mathfrak{A} \left(\mathfrak{U}''(0)\right) \mathfrak{X}(0) + \mathfrak{U}''(0)\right)\\
= & \lambda \mathfrak{X}_{1}(1) + \left(1 - \lambda\right) \mathfrak{X}_{2}(1).
\end{aligned}
\label{Weighted}
\end{equation}

However we note that
$$
\begin{aligned}
& \mathfrak{A}\left(\lambda \mathfrak{U}'(1) + \left(1 - \lambda \right) \mathfrak{U}''(1)\right)\left( \lambda \mathfrak{X}_{1}(1) + \left(1 - \lambda\right) \mathfrak{X}_{2}(1) \right)\\
+ & \lambda\mathfrak{U}'(1)+\left(1 - \lambda\right)\mathfrak{U}''(1) \neq \lambda \mathfrak{X}_{1}(2) + \left(1 - \lambda\right) \mathfrak{X}_{2}(2).
\end{aligned}
$$

Similarly, we have
$$
\begin{aligned}
& \mathfrak{A}\left(\lambda \mathfrak{U}'(k) + \left(1 - \lambda \right) \mathfrak{U}''(k)\right)\left( \lambda \mathfrak{X}_{1}(k) + \left(1 - \lambda\right) \mathfrak{X}_{2}(k) \right)\\
+ & \lambda\mathfrak{U}'(k)+\left(1 - \lambda\right)\mathfrak{U}''(k) \neq \lambda \mathfrak{X}_{1}(k+1) + \left(1 - \lambda\right) \mathfrak{X}_{2}(k+1)
\end{aligned}
$$
for $k > 1$, which completes the proof.
\end{proof}

Lemma \ref{Lemma31} proves that set $\mathcal{U}_{\mathfrak{X}_{T}}$ is not a convex set. Moreover, feasible $\left(\mathfrak{U}(0), \cdots, \mathfrak{U}(K-1)\right) \in \mathcal{U}_{\mathfrak{X}_{T}}$ are solutions of \eqref{momentsystem}, of which the corresponding $u(k), k = 0, \cdots, K-1$ don't have an explicit form of function. Therefore, to obtain an optimal solution to \eqref{optimization} is hardly a possible task.

Due to the complicated topology of the set $\mathcal{U}_{\mathfrak{X}_{T}}$, we don't expect to perform optimization over this set. Instead, we turn our eyes to obtaining a subset of $\mathcal{U}_{\mathfrak{X}_{T}}$ which is convex. By Lemma 2.3 in \cite{wu2023density}, we have that
\begin{equation}
    \e(k_{0}) := \mathfrak{X}_{T} - \mathfrak{X}(k_{0}) \in \mathbb{V}^{2n}_{++}, \quad \exists k_{0} < \infty.
\label{ek}
\end{equation}
Furthermore, we have
\begin{equation}
\mathfrak{X}(k) = \mathfrak{X}(k_{0}) + \omega_{k} \e(k_{0}) \in \mathbb{V}^{2n}_{++}
\label{Xkweight}
\end{equation}
for $k = k_{0}, \cdots, K$ and $0 = \omega_{k_{0}} \leqslant \cdots \leqslant \omega_{K} = 1$. Here the elements of $\mathfrak{X}_{T}$ are the power moments corresponding to the specified terminal distribution $\chi_{T}(\x)$. Next we prove that this choice of $\mathfrak{X}(k)$ essentially minimizes the squared integral distance \cite{hanebeck2014moment} between probability density functions, alongside the system trajectory.

Denote $\alpha_{k} \in \mathbb{R}_{+} \cup \{ +\infty \}, k = k_{0}, \cdots, K-1$ as the weights. We write $\boldsymbol{\alpha}:= \left\{ \alpha_{k_{0}+1}, \cdots, \alpha_{K-1} \right\}$ and $\boldsymbol{\chi}:= \left\{ \chi_{k_{0}+1}, \cdots, \chi_{K-1} \right\}$. Moreover, we let $\alpha_{k_{0}} = 1$. The optimization problem can be formulated as
\begin{equation}
    \min_{\boldsymbol{\chi}}L\left( \boldsymbol{\alpha}, \boldsymbol{\chi} \right)
\label{LalphaChiOpt}
\end{equation}
where
$$
L\left( \boldsymbol{\alpha}, \boldsymbol{\chi} \right) = \sum_{k=k_{0}}^{K-1} \alpha_{k} \int_{\mathbb{R}} \left( \chi_{k+1}(\x) - \chi_{k}(\x) \right)^{2} d\x.
$$

We note that when several $\alpha_{k}$ are $+\infty$, the other terms, of which the weighting coefficients are not these $\alpha_{k}$, are indeed not considered in the cost function.

The directional derivative of $L\left( \boldsymbol{\alpha}, \boldsymbol{\chi} \right)$ with respect to $\delta \chi_{k}$ reads
$$
\begin{aligned}
& \delta L\left( \boldsymbol{\alpha}, \boldsymbol{\chi}; \delta \chi_{k} \right)\\
= & 2\alpha_{k-1} \int_{\mathbb{R}}\delta \chi_{k}\left( \chi_{k}(\x) - \chi_{k-1}(\x) \right)d\x\\
- & 2\alpha_{k} \int_{\mathbb{R}}\delta \chi_{k}\left( \chi_{k+1}(\x) - \chi_{k}(\x) \right)d\x,
\end{aligned}
$$
for $k_{0}+1 \leqslant k \leqslant K-1$. They have to be zero at a minimum for all variations $\delta \chi_{k}$ for $k_{0}+1 \leqslant k \leqslant K-1$. Hence we have
$$
\chi_{k+1} - \chi_{k} = \frac{\alpha_{k-1}}{\alpha_{k}}\left(\chi_{k} - \chi_{k-1}\right)
$$
for $k_{0}+1 \leqslant k \leqslant K-1$. We could then prove that $\chi_{k}$ can be written as 
\begin{equation}
\chi_{k} = \chi_{k_{0}} + \omega_{k}\left( \chi_{K} - \chi_{k_{0}} \right)
\label{ProbPortion}
\end{equation}
for $k_{0}+1 \leqslant k \leqslant K$, where $\omega_{k}$ 
is a function of $\boldsymbol{\alpha}$, namely
$$
\omega_{k} = \frac{1 + \sum_{i=k_{0}+1}^{k}\prod_{j=k_{0}+1}^{i}\frac{\alpha_{j-1}}{\alpha_{j}}}{1 + \sum_{i=k_{0}+1}^{K-1}\prod_{j=k_{0}+1}^{i}\frac{\alpha_{j-1}}{\alpha_{j}}}.
$$
One could also prove that there always exists an $\boldsymbol{\alpha}$ for any arbitrary pair of $\omega_{k_{0}+1:K-1}$ satisfying $0 
\leqslant \omega_{k_{0}+1} \leqslant \cdots \leqslant \omega_{K} = 1$. We note that
$$
\frac{\omega_{k_{0}+2}}{\omega_{k_{0}+1}} = \frac{1+\frac{1}{\alpha_{k_{0}+1}}+\frac{1}{\alpha_{k_{0}+2}}}{1+\frac{1}{\alpha_{k_{0}+1}}} = 1 + \frac{1}{(1+\frac{1}{\alpha_{k_{0}+1}})\alpha_{k_{0}+2}} \geqslant 1.
$$
For any $0 \leqslant \omega_{k_{0}+1} \leqslant \omega_{k_{0}+2} \leqslant 1$, it is obvious that there exists $\alpha_{k_{0}+1}, \alpha_{k_{0}+2}$ satisfying the equality.

Moreover, we have
\begin{equation}
\frac{\omega_{k_{0}+3}}{\omega_{k_{0}+2}} = 1+\frac{1}{(1+\frac{1}{\alpha_{k_{0}+1}}+\frac{1}{\alpha_{k_{0}+2}})\alpha_{k_{0}+3}} \geqslant 1.
\label{k3overk2}
\end{equation}
With $\omega_{k_{0}+1}$ and $\omega_{k_{0}+2}$ given, $(1+\frac{1}{\alpha_{k_{0}+1}}+\frac{1}{\alpha_{k_{0}+2}})$ is determined. It is obvious there exists a positive $\alpha_{k_{0}+3}$ satisfying \eqref{k3overk2} with any $0 \leqslant \omega_{k_{0}+2} \leqslant \omega_{k_{0}+3} \leqslant 1$, given a determined $\omega_{k_{0}+2}$. Similarly, we can prove that there exists a positive $\alpha_{k+1}$ with any $0 \leqslant \omega_{k} \leqslant \omega_{k+1} \leqslant 1$ and $\omega_{k}$ determined, for $k = k_{0}+3, \cdots, K-2$. Therefore, for any $\omega_{k_{0}+1:K-1}$ satisfying $0 \leqslant \omega_{k_{0}+1} \leqslant \cdots \leqslant \omega_{K} = 1$, there exists a corresponding $\boldsymbol{\alpha}$ with all its elements being positive.

For the moment system, from \eqref{ProbPortion}, we have
$$
\int_{\mathbb{R}}\x^{l}\chi_{k}d\x = \int_{\mathbb{R}}\x^{l}\chi_{k_{0}}d\x + \omega_{k}\int_{\mathbb{R}}\x^{l}\left( \chi_{K} - \chi_{k_{0}} \right)d\x
$$
for $l = 1, \cdots, 2n$. Therefore, we have
$$
\mathfrak{X}(k) = \mathfrak{X}(k_{0}) + \omega_{k} \e(k_{0}).
$$
We have proved that the system states of the moment system, in the form of \eqref{Xkweight}, is an optimal solution minimizing a weighted sum of integrated squared distances. By doing this, the variation of distribution of the system state, in the sense of an expected value of squared $L_{2}$ norm, is minimized. In the following parts of this paper, we don't choose $\omega_{k_{0}+1:K-1}$ priorly. Instead, we propose to optimize them according to the criteria for designing the controller.

These results provides us with a way of choosing the subset of $\mathcal{U}_{\mathfrak{X}_{T}}$. Instead of optimizing over all feasible $\mathfrak{U}(k)$, the problem can now be formulated as an optimization over $\omega_{k}$ for $k = k_{0} + 1, \cdots, K - 1$. The advantage of doing this is also obvious: the realizability of $\mathfrak{X}(k)$ for $k = k_{0}, \cdots, K$ is guaranteed, i.e., the Hankel matrices of all $\mathfrak{X}(k)$ are positive-definite. However, the convexity of the set of all feasible $\left(\omega_{k_{0} + 1}, \cdots, \omega_{K-1}\right)$ is not known. Now the task is to prove the convexity of the set of all feasible $\left(\omega_{k_{0}}, \cdots, \omega_{K-1}\right)$.

\begin{proposition}
There exists a sequence 
$$
\left( \breve{\omega}_{k_{0}+1}, \cdots, \breve{\omega}_{K - 1} \right), 0 \leqslant \breve{\omega}_{k} \leqslant 1
$$
for $k = k_{0}, \cdots, K - 1$, with which the following set
$$
\begin{aligned}
& \mathcal{W}_{\mathfrak{X}_{T}} := \\ & \{ \left(\omega_{k_{0} + 1}, \cdots, \omega_{K-1}\right) \mid \omega_{k} \leqslant \breve{\omega}_{k},\\
& \mathfrak{X}(k+1) = \mathfrak{A}(\mathfrak{U}(k))\mathfrak{X}(k)+\mathfrak{U}(k), \\
& \mathfrak{X}(k) = \mathfrak{X}(k_{0}) + \omega_{k - 1} \e(k_{0}), k = k_{0} + 1, \cdots, K - 1,\\
& \omega_{k_{0} + 1} \leqslant \cdots \leqslant \omega_{K-1}\}
\end{aligned}
$$
is convex.
\label{Prop32}
\end{proposition}
\begin{proof}
Substituting $k = K - 1$ and \eqref{Xkweight} into \eqref{momentsystem}, we have
$$
\mathfrak{X}_{T} = \mathfrak{A}(\mathfrak{U}(K-1))\left(\mathfrak{X}_{T} - \left(1 - \omega_{K-1}\right)\e(k_{0})\right)+\mathfrak{U}(K-1).
$$
From Theorem \ref{Theorem31}, by trivially choosing $\omega_{K-1}$ to be zero, $\mathfrak{U}(K - 1) \in \mathbb{V}^{2n}_{++}$. Since $\mathfrak{U}(K - 1)$ is a continuous vector function of $\omega_{K-1}$, there exist subsets of $[0, 1)$ for $\omega_{K-1}$ satisfying that $\left[ \mathfrak{U}(K - 1) \right]_{H} \succ 0$ when $\omega_{K-1}$ falls within these subsets. Next, we further prove that these subsets are path-connected, i.e., the domain of $\omega_{K-1}$ is a convex set.

Assume $\omega'_{K-1}$ and $\omega''_{K-1}$ to fall within these subsets and the corresponding $\mathfrak{U}'(K - 1), \mathfrak{U}''(K - 1) \in \mathbb{V}^{2n}_{++}$. Therefore, there exist two arbitrary $u'(K-1), u''(K-1)$ of which the power moments are $\mathfrak{U}'(K - 1), \mathfrak{U}''(K - 1)$ respectively. Hence by \eqref{uniequation} we can write
$$
\begin{aligned}
x(K)= & a(K-1) \left(\omega'_{K-1}x(k_{0}) + \left(1 - \omega'_{K-1}\right)x(K)\right)\\
+ & u'(K-1)
\end{aligned}
$$
and
$$
\begin{aligned}
x(K)= & a(K-1) \left(\omega''_{K-1}x(k_{0}) + \left(1 - \omega''_{K-1}\right)x(K)\right)\\
+ & u''(K-1).
\end{aligned}
$$
With any arbitrary $\lambda \in \left(0, 1\right)$, we have
$$
\begin{aligned}
x(K) = & \lambda a(K-1) \left(\omega'_{K-1}x(k_{0}) + \left(1 - \omega'_{K-1}\right)x(K)\right)\\
+ & \left(1 - \lambda \right) a(K-1) \left(\omega''_{K-1}x(k_{0}) + \left(1 - \omega''_{K-1}\right)x(K)\right)\\
+ & \lambda u'(K-1) + (1 - \lambda) u''(K-1)\\
= & a(K-1) \left(\left( \lambda \omega'_{K-1} + (1 - \lambda) \omega''_{K-1}\right)x(k_{0})\right.\\
+ & \left. (1 - \lambda \omega'_{K-1} - (1 - \lambda) \omega''_{K-1})x(K)\right)\\
+ & \lambda u'(K-1) + (1 - \lambda) u''(K-1).
\end{aligned}
$$

Therefore, for any $\lambda \omega'_{K-1} + (1 - \lambda) \omega''_{K-1}$ with $\lambda \in [0, 1]$, there always exists a feasible control input $u(k-1) = \lambda u'(K-1) + (1 - \lambda) u''(K-1)$, by which we have $\mathfrak{U}(K - 1) = \lambda \mathfrak{U}'(K-1) + (1 - \lambda) \mathfrak{U}''(K-1) \succ 0$. We then complete the proof that the feasible domain of $\omega_{K-1}$ is convex and can be represented as $\left[0, \breve{\omega}_{K-1} \right)$. Similarly we can prove the feasible domains of $\omega_{k}$ for $k = k_{0}+1, \cdots, K-2$ to be $\left[0, \breve{\omega}_{k} \right)$. In order to ensure the existence of solution to the optimization problem, we consider the closure of the feasible domains, namely $\left[0, \breve{\omega}_{k} \right]$, for $k = k_{0}+1, \cdots, K-1$. Under this circumstance, $\left[\mathfrak{U}(k)\right]_{H}$ is relaxed to a positive semidefinite matrix, namely $\left[\mathfrak{U}(k)\right]_{H} \succeq 0$. By doing this, when $\omega_{k} = \breve{\omega}_{k}$, the distribution of $u(k)$ is discrete, which is supported on finitely many points on $\mathbb{R}$.
\end{proof}

By Proposition \ref{Prop32}, a sub-optimal solution of \eqref{optimization} can then be obtained by the following optimization problem

\begin{equation}
\begin{aligned}
& \operatorname{minimize} f(\omega_{k_{0}}, \cdots, \omega_{K-1})\\
\text{s.t.} \quad & (\omega_{k_{0}}, \cdots, \omega_{K-1}) \in \mathcal{W}_{\mathfrak{X}_{T}}.
\end{aligned}
\label{optimization1}
\end{equation}
It will further become a convex optimization problem, if the function $f(\cdot)$ is chosen as a convex one. In this formulation of the optimization problem, the Hankel matrices of the moment vectors of the system states are confined to be positive definite, which ensures the existence of the system states for the original system \eqref{uniequation}.

\section{Choices of the cost functions}
In the previous section, we proposed an optimization scheme for treating the control of the moment system. However, we have not yet specified the cost function $f(\cdot)$ that we are to use for optimization. In this section, we will put forward different choices of cost functions considering different properties of the control inputs $u(k)$ that we desire, and the corresponding optimization methods.

We note that in the conventional optimal control algorithms, the energy effort is a typical type of cost term, which is the second order moment of a control input. However in our problem, higher order moments are considered for the control task. Different types of cost functions can then be adopted to achieve different design specifications. In the following part of this section, we will propose different design specifications and the corresponding cost functions for the distribution steering problem.

\subsection{Maximal smoothness of state transition}
In our previous paper \cite{wu2023density}, we considered the smoothness of the transition of the system state $\mathfrak{X}(k)$, where we choose $\omega_{k} = \frac{k - k_{0}}{K - k_{0}}$. However, as is mentioned in \cite{wu2023density}, this choice of $\omega_{k}$ doesn't always ensure the positive definiteness of the moment vector $\mathfrak{U}(k)$. We choose the cost function $f$ as
\begin{equation}
f(\omega_{k_{0}}, \cdots, \omega_{K-1}) = \sum_{i = k_{0}}^{K-1}\left( \omega_{i+1} - \omega_{i} \right)^{2} + \omega^{2}_{k_{0}},
\label{costfunc1}
\end{equation}
where $\omega_{K} = 1$. Then we have
$$
\triangledown^{2} f(\omega_{k_{0}}, \cdots, \omega_{K-1}) = \begin{bmatrix}
4 & -2 &  & \\ 
-2 & 4 & -2 & \\ 
 & -2 & \ddots & -2\\ 
 &  & -2 & 4
\end{bmatrix} \succ 0,
$$
i.e., the optimization problem we treat is now a convex one, with the sequence $\left( \omega_{k_{0}}, \cdots, \omega_{K - 1} \right)$ confined to fall within the set $\mathcal{W}_{\mathfrak{X}_{T}}$.

\subsection{Minimum system energy}

In some scenarios, we consider the energy of the system states to be minimized. For example, we consider the cost function
\begin{equation}
f\left(\omega_{k_0}, \cdots, \omega_{K-1}\right)=\sum_{k=k_0}^{K-1} \mathbb{E}\left[x^2(k)\right].
\label{MinSE}
\end{equation}
Now we prove that the RHS of \eqref{MinSE} is also convex. By \eqref{Xkweight}, we have that
$$
\begin{aligned}
& \sum_{k=k_0}^{K-1} \mathbb{E}\left[x^2(k)\right] \\
= & \left(K-k_0\right) \mathbb{E}\left[x^2\left(k_0\right)\right]+\sum_{k=k_0}^{K-1} \omega_k\left(\mathbb{E}\left[x_T^2\right]-\mathbb{E}\left[x^2\left(k_0\right)\right]\right)
\end{aligned}
$$
Since $\mathbb{E}\left[x_T^2\right]$ and $\mathbb{E}\left[x^2\left(k_0\right)\right]$ are constants, we have
$$
\frac{\partial \sum_{k=k_0}^{K-1} \mathbb{E}\left[x^2(k)\right]}{\partial \omega_i \partial \omega_j}=0, \forall k_{0} + 1 \leqslant i, j \leqslant K - 1.
$$
The Hessian matrix is a zero matrix, hence $\sum_{k=k_0}^{K-1} \mathbb{E}\left[x^2(k)\right]$ is convex. Therefore, \eqref{optimization1} with cost function \eqref{MinSE} is a convex optimization problem.

\subsection{Minimum control energy effort}
In some situations, the energy is restricted and we need to take the energy effort into consideration for the control tasks. The cost function can then be chosen as
\begin{equation}
f(\omega_{k_{0}}, \cdots, \omega_{K-1}) = \sum_{i = k_{0}}^{K-1} \mathbb{E}\left[ u^{2}(k) \right]
\label{costfunc2}
\end{equation}

It is a conventional cost function for optimal control. However in our problem, the parameters to be optimized are $\omega_{k}, k = k_{0}, \cdots, K-1$, of which \eqref{costfunc2} is an implicit function. Now it suffices to prove the convexity of \eqref{costfunc2} over $\omega_{k}$.

We first consider $k = k_{0}$. By \eqref{MomentCal} we have
$$
\begin{aligned}
& \mathbb{E}\left[ x(k_{0}+1) \right]\\
= & \mathbb{E}\left[ x(k_{0}) \right] + \omega_{k_{0}}\left( \mathbb{E}\left[ x_{T} \right] - \mathbb{E}\left[ x(k_{0}) \right] \right)\\
= & a(k_{0})\mathbb{E}\left[ x(k_{0}) \right] + \mathbb{E}\left[ u(k_{0}) \right].
\end{aligned}
$$

Then we have
$$
\frac{\partial \mathbb{E}\left[ u(k_{0}) \right]}{\partial \omega_{k_{0}}} = \mathbb{E}\left[ x_{T} \right] - \mathbb{E}\left[ x(k_{0}) \right].
$$

By \eqref{MomentCal} we could also write
$$
\begin{aligned}
& \mathbb{E}\left[ x^{2}(k_{0}+1) \right]\\
= & \mathbb{E}\left[ x^{2}(k_{0}) \right] + \omega_{k_{0}}\left( \mathbb{E}\left[ x^{2}_{T} \right] - \mathbb{E}\left[ x^{2}(k_{0}) \right] \right)\\
= & a^{2}(k_{0})\mathbb{E}\left[ x^{2}(k_{0}) \right] + 2a(k_{0})\mathbb{E}\left[ x(k_{0}) \right]\mathbb{E}\left[ u(k_{0}) \right]\\
+ & \mathbb{E}\left[ u^{2}(k_{0}) \right].
\end{aligned}
$$

Now the second order moment of $u(k_{0})$ reads
$$
\begin{aligned}
& \mathbb{E}\left[ u^{2}(k_{0}) \right]\\
= & \left( 1 - a^{2}\left( k_{0} \right) \right)\mathbb{E}\left[ x^{2}(k_{0}) \right] - 2a(k_{0})\mathbb{E}\left[ x(k_{0}) \right]\mathbb{E}\left[ u(k_{0}) \right]\\
+ & \omega_{k_{0}}\left( \mathbb{E}\left[ x^{2}_{T} \right] - \mathbb{E}\left[ x^{2}(k_{0}) \right] \right).
\end{aligned}
$$

By differentiating both sides of the equation over $\omega_{k_{0}}$, we have
$$
\begin{aligned}
& \frac{\partial \mathbb{E}\left[ u^{2}(k_{0}) \right]}{\partial \omega_{k_{0}}}\\
= & - 2a(k_{0})\mathbb{E}\left[ x(k_{0}) \right]\frac{\partial\mathbb{E}\left[  u(k_{0}) \right]}{\partial \omega_{k_{0}}}\\
+ & \mathbb{E}\left[ x^{2}_{T} \right] - \mathbb{E}\left[ x^{2}(k_{0}) \right]\\
= & - 2a(k_{0})\mathbb{E}\left[ x(k_{0}) \right]\left(\mathbb{E}\left[ x_{T} \right] - \mathbb{E}\left[ x(k_{0}) \right]\right)\\
+ & \mathbb{E}\left[ x^{2}_{T} \right] - \mathbb{E}\left[ x^{2}(k_{0}) \right].
\end{aligned}
$$
Since $a(k_{0}), \mathbb{E}\left[ x(k_{0}) \right], \mathbb{E}\left[ x_{T} \right], \mathbb{E}\left[ x^{2}_{T} \right], \mathbb{E}\left[ x^{2}(k_{0}) \right]$ are all constant scalars or vectors, we have
\begin{equation}
\frac{\partial \mathbb{E}\left[ u^{2}(k_{0}) \right]}{\partial \omega^{2}_{k_{0}}} = 0.
\label{k0partialpos}
\end{equation}

Then we consider $k_{0} + 1 \leqslant k \leqslant K - 1$. Since
$$
\begin{aligned}
& \mathbb{E}\left[ x(k+1) \right]\\
= & \mathbb{E}\left[ x(k_{0}) \right] + \omega_{k}\left( \mathbb{E}\left[ x_{T} \right] - \mathbb{E}\left[ x(k_{0}) \right] \right)\\
= & a(k)\mathbb{E}\left[ x(k) \right] + \mathbb{E}\left[ u(k) \right]\\
= & a(k)\left( \mathbb{E}\left[ x(k_{0}) \right] + \omega_{k-1}\left( \mathbb{E}\left[ x_{T} \right] - \mathbb{E}\left[ x(k_{0}) \right] \right) \right) + \mathbb{E}\left[ u(k) \right],
\end{aligned}
$$
which yields
$$
\frac{\partial \mathbb{E}\left[ u(k) \right]}{\partial \omega_{k}} = \mathbb{E}\left[ x_{T} \right] - \mathbb{E}\left[ x(k_{0}) \right]
$$
and
$$
\frac{\partial \mathbb{E}\left[ u(k) \right]}{\partial \omega_{k-1}} = -a(k) \left(\mathbb{E}\left[ x_{T} \right] - \mathbb{E}\left[ x(k_{0}) \right]\right).
$$

The second order moment of $x(k+1)$ reads
$$
\begin{aligned}
& \mathbb{E}\left[ x^{2}(k+1) \right]\\
= & \mathbb{E}\left[ x^{2}(k_{0}) \right] + \omega_{k}\left( \mathbb{E}\left[ x^{2}_{T} \right] - \mathbb{E}\left[ x^{2}(k_{0}) \right] \right)\\
= & a^{2}(k)\mathbb{E}\left[ x^{2}(k_{0}) \right] + a^{2}(k)\omega_{k-1}\left( \mathbb{E}\left[ x^{2}_{T} \right] - \mathbb{E}\left[ x^{2}(k_{0}) \right] \right)\\
+ & 2a(k)\omega_{k-1}\left( \mathbb{E}\left[ x_{T} \right] - \mathbb{E}\left[ x(k_{0}) \right] \right)\mathbb{E}\left[ u(k) \right]\\
+ & 2a(k)\mathbb{E}\left[ x(k_{0}) \right]\mathbb{E}\left[ u(k) \right] + \mathbb{E}\left[ u^{2}(k) \right].
\end{aligned}
$$

It is easy to verify that
\begin{equation}
\begin{aligned}
& \frac{\partial \mathbb{E}\left[ u^{2}(k) \right]}{\partial \omega_{i} \partial \omega_{j}}\\
= & \left\{\begin{matrix}
2a^{2}(k)\left( \mathbb{E}\left[ x^{2}_{T} \right] - \mathbb{E}\left[ x^{2}(k_{0}) \right] \right)^{2} & i=j=k\\ 
-2a(k)\left( \mathbb{E}\left[ x^{2}_{T} \right] - \mathbb{E}\left[ x^{2}(k_{0}) \right] \right)^{2} & i=k, j=k+1\\ 
0 & \text{otherwise}.
\end{matrix}\right.
\end{aligned}
\end{equation}

We note that the elements of the Hessian of $f$ in \eqref{costfunc2} are functions of $a(k)$ for $k = k_{0}, \cdots, K-1$, which makes the positive semidefiniteness of the Hessian not ensured. The convexity of $f$ is then not guaranteed either. On the other hand, we proved $\mathcal{W}_{\mathfrak{X}_{T}}$ is convex and compact in Proposition \ref{Prop32}. It then makes it feasible to discretize $\mathcal{W}_{\mathfrak{X}_{T}}$ into grids and use a exhaustive search method to treat the optimization problem \eqref{optimization1}.

\subsection{A more general form of cost function}

We consider a more general form of cost function. The cost function reads
\begin{equation}
\begin{aligned}
& f\left(\omega_{k_0}, \cdots, \omega_{K-1}\right) \\
= & \mathbb{E}\left[\alpha \mathcal{K}\left(x(k)\right)+\beta \mathcal{G}\left(u(k)\right)\right]\\
= & \alpha \mathbb{E}\left[\mathcal{K}\left(x(k)\right)\right]+\beta \mathbb{E}\left[\mathcal{G}\left(u(k)\right)\right]\\
\end{aligned}
\label{costfunc3}
\end{equation}
where $\alpha, \beta > 0$ are weights, and $\mathcal{K}(\cdot), \mathcal{G}(\cdot)$ are polynomial functions. For this type of cost function, an exhaustive search method, as was used in Subsection C, can be applied to the optimization problem.

In this paper, we mainly consider the previous four cost functions. However, the cost functions are not limited to these four. Cost functions considering other orders of power moments can also be applied to form the optimization problem.

\section{Realization of the control inputs and an algorithm for distribution steering}

In the previous section, we put forward a control law for the moment system in the manner of the conventional optimal control scheme. However by the control law in the previous sections, the control inputs we obtained are those of the moment system, i.e., $\mathfrak{U}(k)$ for $k=0, \cdots, K-1$. In order to control the primal system \eqref{uniequation}, we need to further obtain $u(k)$ for $k=0, \cdots, K-1$. In this section, we will propose an algorithm to determine the $u(k)$ given $\mathfrak{U}(k)$ obtained by the optimization problem (22). This problem is an ill-posed one, i.e., there might be infinitely many feasible $u(k)$ for a given $\mathfrak{U}(k)$. However, we will select a unique solution $u(k)$ by the algorithm proposed in this section, which satisfies the given $\mathfrak{U}(k)$. That's why we use the word determine here.

Moreover, for the sake of simplicity, we omit $k$ if there is no ambiguity in the following part of this section. The problem now becomes that of proposing an algorithm which estimates the distribution of $u(k)$, for which the power moments are as specified.

A convex optimization scheme for distribution estimation by the Kullback-Leibler distance has been proposed in \cite{wu2023non} considering the Hamburger moment problem, which is used for control input realization in our previous paper \cite{wu2023density}. Moreover, we observed that the performance of estimation for probability distributions which are relatively smooth can be improved by using the squared Hellinger distance as the metric \cite{wu2022non1}. We adopt this strategy in this paper for treating the realization of the control inputs. The procedure is as follows.
Let $\mathcal{P}$ be the space of probability distributions on the real line with support there, and let $\mathcal{P}_{2n}$ be the subset of all $p \in \mathcal{P}$ which have at least $2 n$ finite moments (in addition to $\mathbb{E}\left[u^0(k)\right]$, which of course is 1). The squared Hellinger distance is then defined as
\begin{equation}
\mathbb{H}^2(\theta, \nu)=\int_{\mathbb{R}}(\sqrt{\theta(\u)}-\sqrt{\nu(\u)})^2 d \u
\label{SHellinger}
\end{equation}
where $\theta$ is an arbitrary probability distribution in $\mathcal{P}$. Moreover, we define
$$
K(\u)=\left[\begin{array}{lllll}
1 & \u & \cdots & \u^{n-1} & \u^n
\end{array}\right]^T
$$
and
$$
\boldsymbol{\mathscr{M}}_{2n}=\left[\begin{array}{cccc}
1 & \mathbb{E}[u] & \cdots & \mathbb{E}\left[u^n\right] \\
\mathbb{E}[u] & \mathbb{E}\left[u^2\right] & \cdots & \mathbb{E}\left[u^{n+1}\right] \\
\vdots & \vdots & \ddots & \\
\mathbb{E}\left[u^n\right] & \mathbb{E}\left[u^{n+1}\right] & & \mathbb{E}\left[u^{2 n}\right]
\end{array}\right]
$$
where $\mathbb{E}\left[u^i\right], i=1, \cdots, 2 n$ are the elements of the control vector $\mathfrak{U}$ of the moment system.

We define the linear integral operator $\Xi$ as
$$
\Xi: \nu(\u) \mapsto \boldsymbol{\mathscr{M}}_{2n}=\int_{\mathbb{R}} K(\u) \nu(\u) K^T(\u) d \u,
$$
where $\nu(\u)$ belongs to the space $\mathcal{P}_{2n}$. Furthermore, the range $(\Xi)=\Xi \mathcal{P}_{2n}$ is also convex since $\mathcal{P}_{2n}$ is convex. We let
$$
\mathcal{S}_{+}:=\left\{\Lambda \in \operatorname{range}(\Xi) \mid K(\u)^T \Lambda K(\u)>0, \u \in \mathbb{R}\right\}.
$$

Given any $\theta \in \mathcal{P}$ and any $\boldsymbol{\mathscr{M}}_{2n} \succ 0$, there is a unique $\hat{\nu} \in \mathcal{P}_{2n}$ that minimizes \eqref{SHellinger} subject to $\Xi(\hat{\nu})=\boldsymbol{\mathscr{M}}_{2n}$, namely
$$
\hat{\nu}=\frac{\theta}{\left(1+K^T \hat{\Lambda} K\right)^2}
$$
where $\hat{\Lambda}$ is the unique solution to the problem of minimizing
\begin{equation}
\mathbb{J}_{\theta}(\Lambda):=\operatorname{tr}(\Lambda \boldsymbol{\mathscr{M}}_{2n})+\int_{\mathbb{R}} \frac{\theta}{1+K^T \Lambda K} d \u.
\label{Jr}
\end{equation}

Then the distribution estimation is formulated as a convex optimization problem. The map $\Lambda \mapsto \boldsymbol{\mathscr{M}}_{2n}$ is proved to be homeomorphic, which ensures the existence and uniqueness of the solution to the realization of control inputs \cite{wu2022non1}. Unlike other moment methods, the power moments of our proposed distribution estimate are exactly identical to those specified, which makes it a satisfactory approach for realization of the control inputs \cite{wu2022non1}. Since the prior distribution $\theta(\u)$ and the distribution estimate $\nu(\u)$ are both supported on $\mathbb{R}$, $\theta(\u)$ can be chosen as a Gaussian distribution (or a Cauchy distribution if $\hat{\nu}(\u)$ is assumed to be heavy-tailed). 

\section{Two types of general distribution steering problems and the corresponding algorithms}

In the previous sections of the paper, we considered the general distribution steering problem which only assumes the existence of the first several finite power moments. Loosely speaking, the distributions can be divided into two types, namely the continuous and discrete ones. In this section, we will propose algorithms corresponding to the two types of distributions.

\subsection{An algorithm for continuous distribution steering}

We first consider the continuous distribution steering algorithm, which is concluded in the following Algorithm \ref{alg:1}.

\begin{algorithm}[htbp]
    \caption{Continuous distribution steering.}
    \label{alg:1}
    \begin{algorithmic}[1]
        \Require The maximal time step $K$; the parameter of the system $a(k)$ for $k = 0, \cdots, K-1$; the initial system distribution $\chi_{0}(\x)$; the specified terminal distribution $\chi_{T}(\x)$.
        \Ensure The controls $u(k)$, $k = 0, \cdots, K-1$.
        \State $k \Leftarrow 0$
    \While{$k < K$ and $\e(k) \notin \mathbb{V}^{2n}_{++}$}
        \State Calculate $\mathfrak{X}(k)$ by \eqref{momentsystem} if $k > 0$ or by \eqref{XK} if $k = 0$
        \State Calculate $\e(k)$ by \eqref{ek}
        \If{$\e(k) \in \mathbb{V}^{2n}_{++}$}
        \State Optimize the cost function $f\left(\omega_{k_0}, \cdots, \omega_{K-1}\right)$ over the domain $\mathcal{W}_{\mathfrak{X}_T}$. Obtain the optimal $\omega_{k_0}^*, \cdots, \omega_{K-1}^*$.
        \State Calculate the states of the moment system $\mathfrak{X}(i)$ for $i = k + 1, \cdots, K-1$ by \eqref{Xkweight} with $\omega_{k_0}^*, \cdots, \omega_{K-1}^*$ \label{Step6}
        \State Calculate the controls of the moment system $\mathfrak{U}(i)$ for $i = k, \cdots, K-1$ by \eqref{momentsystem}
        \State Optimize the cost function \eqref{Jr} and obtain the analytic estimates of the distributions $\hat{\nu}_{i}(\u)$ for $i = k, \cdots, K-1$
        \Else
            \State $u(k) = 0$
        \EndIf
        \State Calculate the power moments of the system state $x(k+1)$, i.e., $\mathfrak{X}(k+1)$
        \State $k \Leftarrow k+1$
    \EndWhile
    \end{algorithmic}
\end{algorithm}

There is still an important issue to consider in the algorithm, which is to determine the set $\mathcal{W}_{\mathfrak{X}_T}$. By the proof of Proposition \ref{Prop32}, it is equivalent to determine the maximal $\omega_{K-1} \in(0,1)$. It can be treated by the following optimization.
$$
\begin{array}{ll} 
& \max \omega_{K-1} \\
\text { s.t. } & \mathfrak{X}_T=\mathfrak{A}(\mathfrak{U}(K-1))\left(1-\omega_{K-1}\right) \mathfrak{X}_T+\mathfrak{U}(K-1), \\
& {\left[\begin{array}{ccc}
1 & \cdots & \mathbb{E}\left[u^n(K-1)\right] \\
\vdots & \ddots & \\
\mathbb{E}\left[u^n(K-1)\right] & & \mathbb{E}\left[u^{2 n}(K-1)\right]
\end{array}\right] \succeq 0} \\
& 0<\omega_{K-1}<1 .
\end{array}
$$

We take $\breve{\omega}_{K-1} = \max \omega_{K-1}$, and obtain the corresponding $\breve{\mathfrak{X}}(K-1)$. Similarly, we obtain $\breve{\omega}_{k} = \max \omega_{k}$ and the corresponding $\breve{\mathfrak{X}}(k)$ by the convex optimization
$$
\begin{array}{ll} 
& \max \omega_{k} \\
\text { s.t. } & \breve{\mathfrak{X}}_{k+1}=\mathfrak{A}(\mathfrak{U}(k))\left(1-\omega_{k}\right) \mathfrak{X}_{k}+\mathfrak{U}(k), \\
& {\left[ \mathfrak{U}(k) \right] \succeq 0} \\
& 0<\omega_{k}<\breve{\omega}_{k+1}
\end{array}
$$
in an reversed order from $k=K-2$ to $k_{0}$.

\subsection{An algorithm for discrete distribution steering}
In the real applications, we are sometimes confronted with the problem of steering a colossal group of discrete agents, which are distributed arbitrarily in the whole domain rather than following a prescribed distribution. Considering this type of problem, we characterize the distribution of the agents as an occupation measure following \cite{wu2022group}

$$
d \chi_{k}(\x, \u)=\frac{1}{N} \sum_{i=1}^{N} \delta\left(\x-x_{i}(k)\right) \delta\left(\u-u_{i}(k)\right) d \x d \u,
$$

Then we can write the power moments of the system states and control inputs respectively as
\begin{equation}
\begin{aligned}
\mathbb{E}\left[ x^{l}(k) \right]
& = \int_{\mathbb{R} \times \mathbb{R}} \x^{l} d \chi_{k}(\x, \u)\\
& = \int_{\mathbb{R}} \x^{l} \sum_{i=1}^{N} \delta\left(\x-x_{i}(k)\right)d\x\\
& = \frac{1}{N}\sum_{i=1}^{N}x^{l}_{i}(k),
\label{xlK1}
\end{aligned}
\end{equation}
and
\begin{equation}
\mathbb{E}\left[ \u^{l}(k) \right] = \int_{\mathbb{R} \times \mathbb{R}} \u^{l} d \chi_{k}(\x, \u) = \frac{1}{N}\sum_{i=1}^{N}u^{l}_{i}(k).
\label{ulK1}
\end{equation}

The occupation measure steering problem differs from the continuous distribution steering one mainly in determining the control inputs for each agent, which means that we have to draw samples from the realized control inputs. Since the realized controls by our proposed algorithm have analytic form of function, acceptance-rejection sampling \cite{casella2004generalized} strategy can be used for this task. The idea of acceptance-rejection sampling is that even it is not feasible for us to directly sample from the functions of the control inputs, there exists another candidate distribution, from which it is easy to sample from. A common choice of light-tailed distributions is the Gaussian. Then the task can be reduced to sampling from the candidate distribution directly and then rejecting the samples in a strategic way to make the remaining samples seemingly drawn from the distributions of the control inputs.

By adopting the acceptance-rejection sampling strategy, we update Algorithm \ref{alg:1} as to treat the occupation measure steering problem, which is given in Algorithm \ref{alg:3}. 

\begin{algorithm}
    \caption{Discrete distribution steering.}
    \label{alg:3}
    \begin{algorithmic}[1]
        \Require The number of agents $N \in \mathbb{N}_{0}$; the maximal time step $K$; the parameter of the system $a(k)$ for $k = 0, \cdots, K-1$; the initial occupation measure $d\chi_{0}(\x)$; the specified terminal occupation measure $d\chi_{T}(\x)$.
        \Ensure The control inputs for the $i_\text{th}$ target $u_{i}(k)$, $k = 0, \cdots, K-1$, $i = 1, \cdots, N$.
        \State $k \Leftarrow 0$
    \While{$k < K$ and $\e(k) \notin \mathbb{V}^{2n}_{++}$}
        \State Calculate $\mathfrak{X}(k)$ by \eqref{momentsystem} if $k > 0$ or by \eqref{XK} if $k = 0$
        \State Calculate $\e(k)$ by \eqref{ek}
        \If{$\e(k) \in \mathbb{V}^{2n}_{++}$}
        \State Optimize the cost function $f\left(\omega_{k_0}, \cdots, \omega_{K-1}\right)$ over the domain $\mathcal{W}_{\mathfrak{X}_T}$. Obtain the optimal $\omega_{k_0}^*, \cdots, \omega_{K-1}^*$
        \State Calculate the states of the moment system $\mathfrak{X}(i)$ for $i = k + 1, \cdots, K-1$ by \eqref{Xkweight} with $\omega_{k_0}^*, \cdots, \omega_{K-1}^*$ \label{Step61}
        \State Calculate the controls of the moment system $\mathfrak{U}(i)$ for $i = k, \cdots, K-1$ by \eqref{momentsystem}
        \State Optimize the cost function \eqref{Jr} and obtain the analytic estimates of the distributions $\hat{\nu}_{i}(\u)$ for $i = k, \cdots, K-1$
        \State Sample the control inputs $u_{i}(j)$ of all agents at time step $j = k, \cdots, K-1$ by the acceptance-rejection strategy.
        \Else
            \State $u_{i}(k) = 0, i = 1, \cdots, N$
        \EndIf
        \State Calculate the power moments of the system state $x(k+1)$, i.e., $\mathfrak{X}(k+1)$
        \State $k \Leftarrow k+1$
    \EndWhile
    \end{algorithmic}
\end{algorithm}

\section{The power moments and the Gaussian mixture model: A comparison}

The general distribution steering problem has recently garnered significant interest in the control community. As an infinite-dimensional problem, addressing it requires dimension reduction. In this paper, we propose using the truncated moment sequence to characterize infinite-dimensional probability distributions. In addition to using moments, Gaussian Mixture Models (GMMs), a widely used algorithm for approximating probability distributions, have also been considered for general distribution steering tasks \cite{balci2024density, kumagai2024chance}.
Although GMMs simplify the treatment of infinite-dimensional distribution problems and can leverage results from existing literature on distribution steering for Gaussian distributions, they do not necessarily converge to the true distribution. Notably, the decay rate of the distribution approximation by GMMs is always squared exponential. This means that even as the number of Gaussian components approaches infinity, the mixture may fail to converge to the true distribution if the decay rate of the true distribution is not squared exponential. General distributions usually have all finite power moments but do not exhibit squared exponential decay. Common examples include the exponential, Gamma, Pareto, and Chi-squared distributions, of which the decay rate is less than the squared exponential. As a result, GMMs perform poorly for these types of probability distributions. Moreover, the GMMs are usually determined by expectation maximization \cite{moon1996expectation}, in which the optimization problems are not convex, which would also effect the performance. In contrast, the algorithm we propose ensures that the power moments of the terminal distribution are identical to those of the desired distribution (the control inputs of the moment system are obtained by convex optimization with solution being proved to exist and are unique). As the number of moment terms increases, specifically as $2n \rightarrow +\infty$, $\chi_{K}$ almost surely converges to the true distribution \cite[Theorem 4.5.5]{chung2001course}. If both the initial and desired terminal distributions are continuous, the terminal distribution converges to the desired one. Furthermore, the error in the terminal distribution, measured in terms of total variation distance, uniformly converges to zero. As emphasized earlier, the general distribution steering problem is infinite-dimensional, and some error between the terminal and desired distributions is inevitable for any finite $n$. In our previous work \cite{wu2022group}, we derived a tight upper bound for this error in terms of total variation distance, which remains valid for the realization of control inputs using the squared Hellinger distance in this paper. This property distinguishes our proposed algorithm from the GMM-based approach.

\section{Numerical results and comparison between cost functions}

In this section, we will simulate general distribution steering problems with the cost functions proposed in the previous sections of the paper. We consider two typical scenarios in real applications. The first one is to separate a group of agents into several smaller groups. The second one is to steer the agents which are in separate groups to desired terminal groups. For the first type of problem, we consider to steer a Gaussian distribution to a mixture of two generalized logistic distributions as an example. And for the second type of problem, we consider to steer a mixture of two Laplacians to a mixture of two Gaussians.

\subsection{A Gaussian to two generalized logistic distributions}
We first consider the problem of steering a Gaussian distribution to a mixture of generalized logistic distributions with two modes. The initial one is chosen as
\begin{equation}
\chi_0(\x)=\frac{1}{\sqrt{2 \pi}} e^{\frac{x^2}{2}}
\label{q01}
\end{equation}
and the terminal one is specified as
\begin{equation}
\chi_T(\x)=\frac{0.3 \cdot 2e^{-x+2}}{(1 + e^{-x+2})^{3}} + \frac{0.7 \cdot 3e^{-x - 3}}{(1 + e^{-x - 3})^{4}}.
\label{qt1}
\end{equation}
The system parameters $a(k), k=0, \cdots, 3$ are i.i.d. samples drawn from the uniform distribution $U[0.3,0.5]$. The dimension of each $\mathfrak{U}(k)$ is $4$. 

\begin{figure}[htbp]
\centering
\includegraphics[scale=0.30]{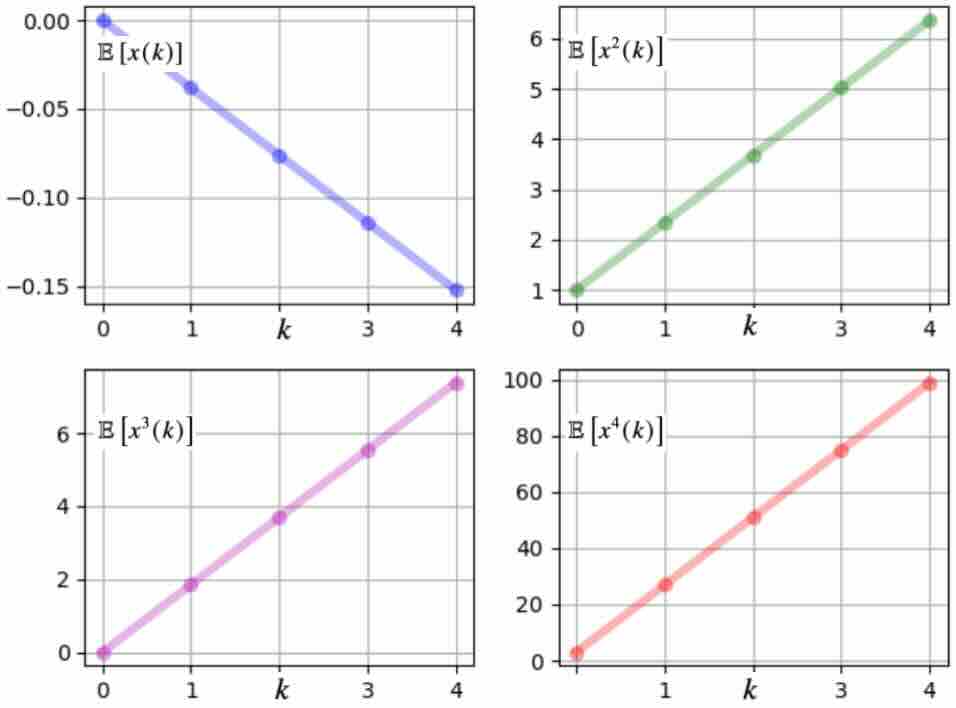}
\centering
\caption{$\mathfrak{X}(k)$ at time steps $k = 0, 1, 2, 3, 4$ with cost function \eqref{costfunc1}. The upper left figure shows $\mathbb{E}\left[ x(k)\right]$. The upper right one shows $\mathbb{E}\left[ x^{2}(k)\right]$. The lower left one shows $\mathbb{E}\left[ x^{3}(k)\right]$ and the lower right one shows $\mathbb{E}\left[ x^{4}(k)\right]$.}
\label{fig1}
\end{figure}

\begin{figure}[htbp]
\centering
\includegraphics[scale=0.30]{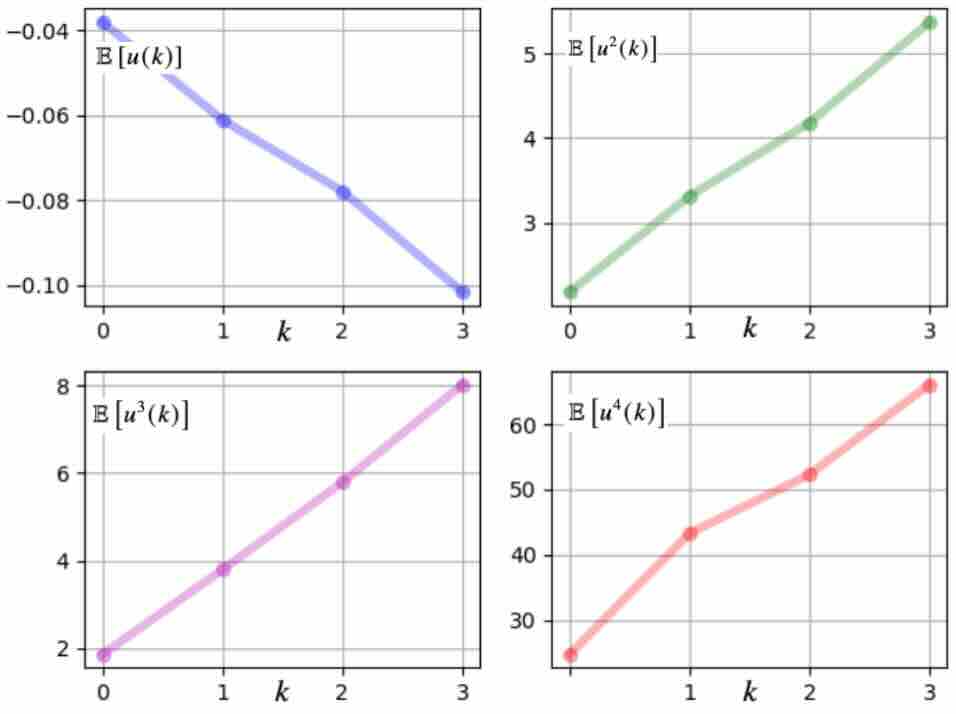}
\centering
\caption{$\mathfrak{U}(k)$ at time steps $k = 0, 1, 2, 3$ with cost function \eqref{costfunc1}. The upper left figure shows $\mathbb{E}\left[ u(k)\right]$. The upper right one shows $\mathbb{E}\left[ u^{2}(k)\right]$. The lower left one shows $\mathbb{E}\left[ u^{3}(k)\right]$ and the lower right one shows $\mathbb{E}\left[ u^{4}(k)\right]$.}
\label{fig2}
\end{figure}

\begin{figure}[htbp]
\centering
\includegraphics[scale=0.30]{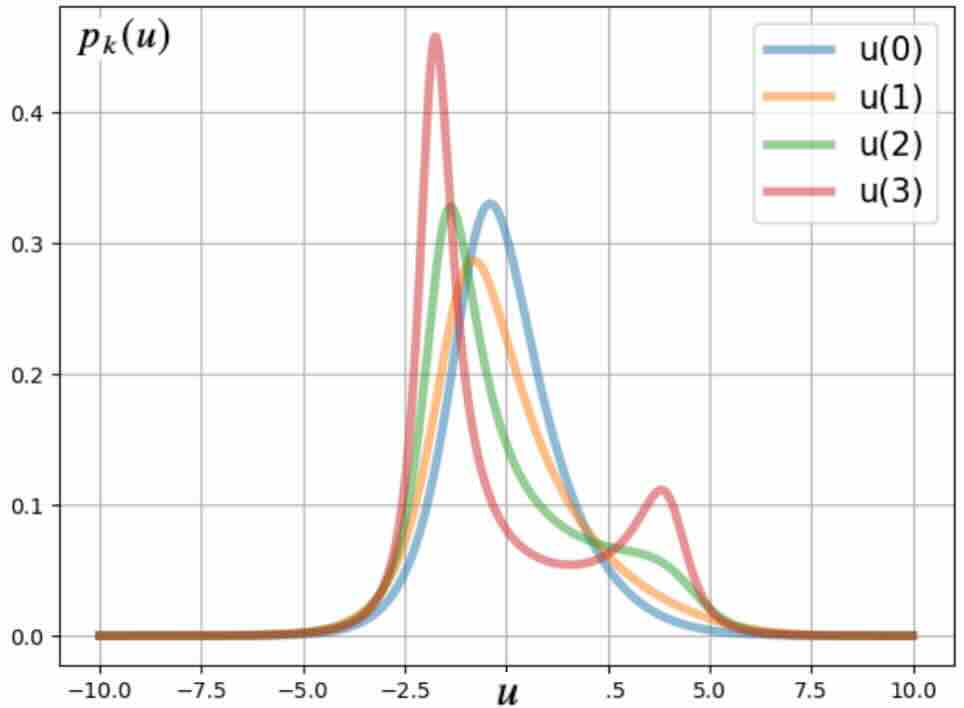}
\centering
\caption{Realized control inputs $\nu_{k}(\u)$ of $u(k)$ by $\mathfrak{U}(k)$ for $k = 0, 1, 2, 3$ , which are obtained by cost function \eqref{costfunc1}.}
\label{fig3}
\end{figure}

We first consider the maximal smoothness of state transition as the control criterion, i.e., choose the cost function as \eqref{costfunc1}. The states of the moment system, i.e., $\mathscr{X}(k)$ for $k=0,1,2,3$, are illustrated in Figure \ref{fig1}. The controls of the moment system, i.e., $\mathfrak{U}(k)$ for $k=0,1,2,3$ are given in Figure \ref{fig2}. The realized controls in Figure \ref{fig3} also show that the transition of the control inputs is smooth, even the specified terminal distribution has two modes, which are Laplacians. However, the tradeoff of the smooth transition is a relatively large energy effort $\sum_{k=0}^{3}\mathbb{E}\left[ u^{2}(k) \right] = 14.988$.

\begin{figure}[htbp]
\centering
\includegraphics[scale=0.30]{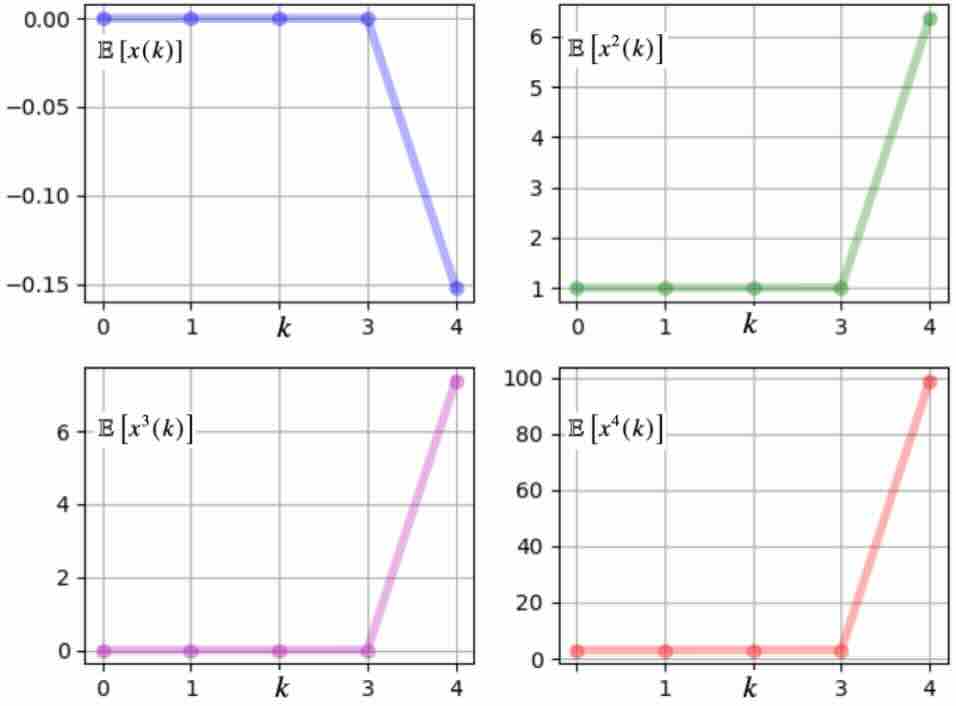}
\centering
\caption{$\mathfrak{X}(k)$ at time steps $k = 0, 1, 2, 3, 4$ with cost function \eqref{costfunc2}.}
\label{fig4}
\end{figure}

\begin{figure}[htbp]
\centering
\includegraphics[scale=0.30]{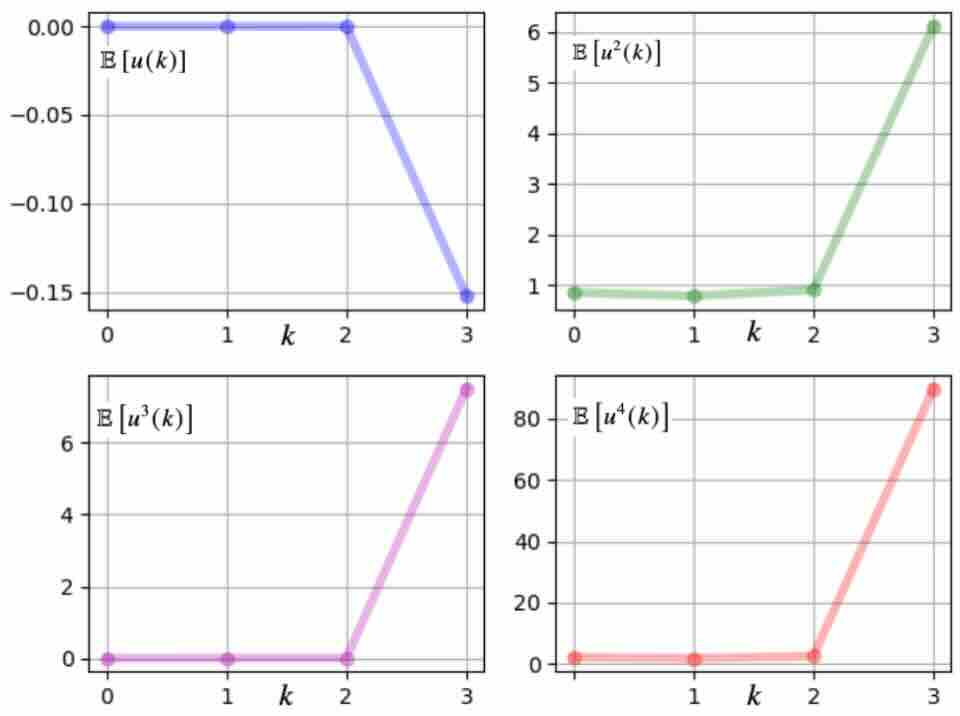}
\centering
\caption{$\mathfrak{U}(k)$ at time steps $k = 0, 1, 2, 3$ with cost function \eqref{costfunc2}.}
\label{fig5}
\end{figure}

\begin{figure}[htbp]
\centering
\includegraphics[scale=0.30]{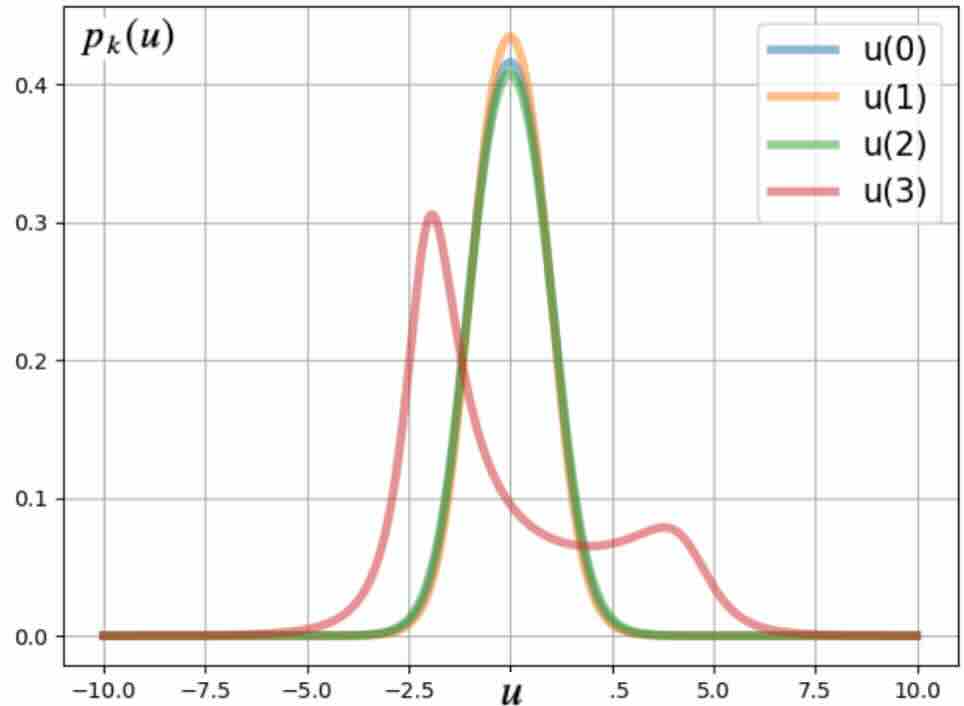}
\centering
\caption{Realized control inputs $\nu_{k}(\u)$ of $u(k)$ by $\mathfrak{U}(k)$ for $k = 0, 1, 2, 3$, which are obtained by cost function \eqref{costfunc2}.}
\label{fig6}
\end{figure}

In particular circumstances, the energy effort we are able to provide is quite limited. For the distribution steering problems which are sensitive to energy, we choose the cost function as \eqref{costfunc2}. The optimization is then performed by discretizing the whole domain and searching exhaustively on it. The size of grid of each dimension $\omega_{k}$ is $0.01$. The results are given in Figure \ref{fig4}, \ref{fig5} and \ref{fig6}. We note that the transition of the system state is not quite smooth as shown in Figure \ref{fig6}. However, the energy effort $\sum_{k=0}^{3}\mathbb{E}\left[ u^{2}(k) \right] = 9.007$, which is much less than that by using the smoothness of state transition as the cost function. 

In situations where both smoothness of the control inputs and the energy effort are considered, the cost function \eqref{costfunc3} provides us with a treatment to the distribution steering problem. In this simulation, we choose the cost function as
\begin{equation}
\begin{aligned}
& f\left(\omega_{0}, \cdots, \omega_{3}\right)\\
= & \mathbb{E}\left[u^{2}(0) \right] + \mathbb{E}\left[u^{2}(1) \right] + 4 \mathbb{E}\left[u^{2}(2) \right]\\
+ & 18 \mathbb{E}\left[u^{2}(3) \right] + \frac{2}{5} \sum_{k=0}^{3}\mathbb{E}\left[ x^{2}(k) \right].
\end{aligned}
\label{WeightedCost}
\end{equation}

\begin{figure}[htbp]
\centering
\includegraphics[scale=0.30]{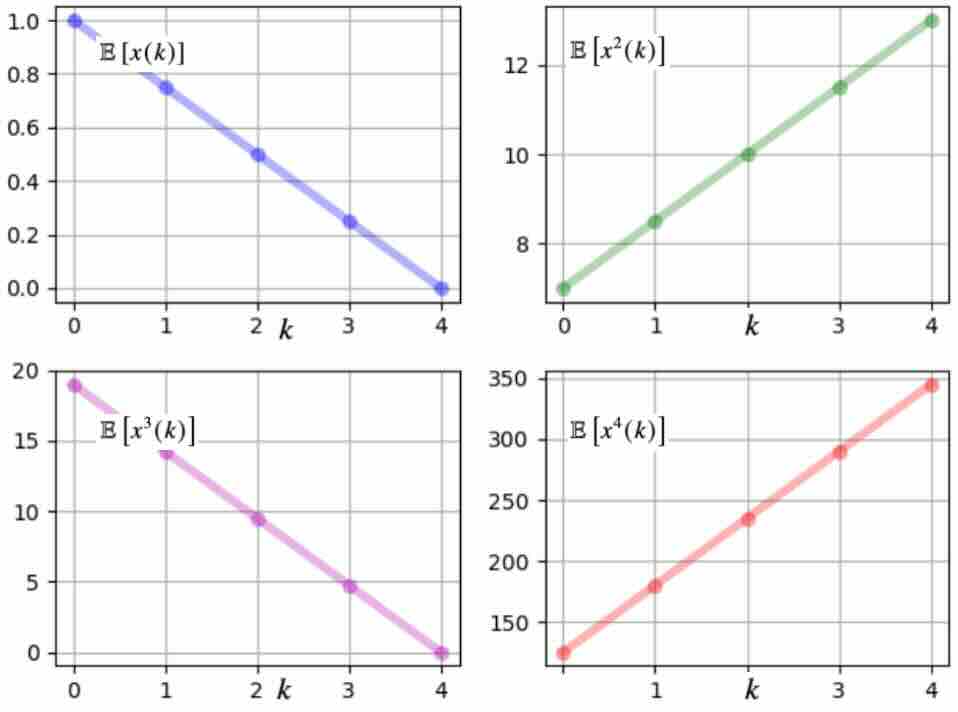}
\centering
\caption{$\mathfrak{X}(k)$ at time steps $k = 0, 1, 2, 3, 4$ with cost function \eqref{WeightedCost}.}
\label{fig7}
\end{figure}

\begin{figure}[htbp]
\centering
\includegraphics[scale=0.30]{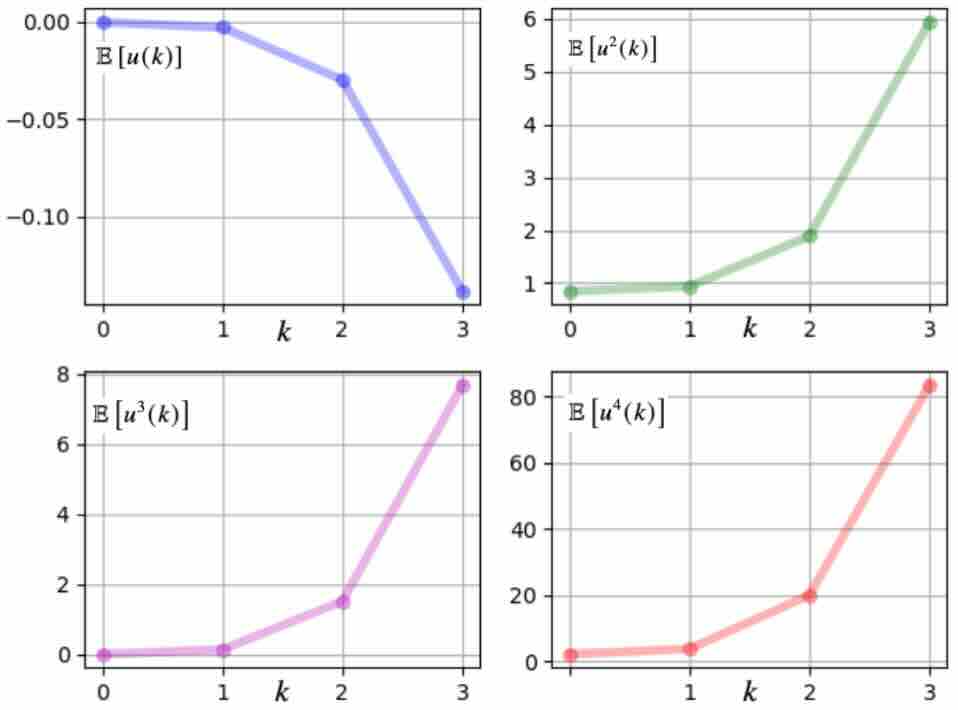}
\centering
\caption{$\mathfrak{U}(k)$ at time steps $k = 0, 1, 2, 3$ with cost function \eqref{WeightedCost}.}
\label{fig8}
\end{figure}

\begin{figure}[htbp]
\centering
\includegraphics[scale=0.30]{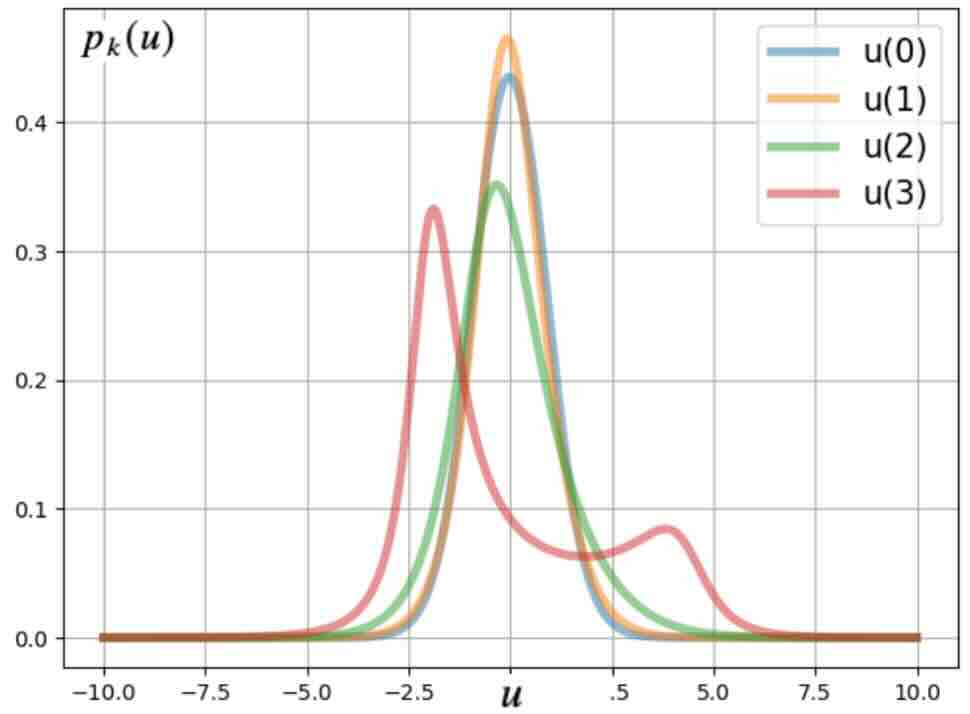}
\centering
\caption{Realized control inputs $\nu_{k}(\u)$ of $u(k)$ by $\mathfrak{U}(k)$ for $k = 0, 1, 2, 3$, which are obtained by cost function \eqref{WeightedCost}.}
\label{fig9}
\end{figure}

The size of grid of each dimension is $0.01$. The simulation results are given in Figure \ref{fig7}, \ref{fig8} and \ref{fig9}. We note that the transition of the control inputs are smoother than the distribution steering by merely considering the energy effort. The energy effort $\sum_{k=0}^{3}\mathbb{E}\left[ u^{2}(k) \right] = 10.648$, which is larger compared to that obtained by \eqref{costfunc2} however is relatively smaller than that obtained by \eqref{costfunc1}. The cost function, in the form of a weighted mixture of the energy effort and the system energy, provides us with a balanced choice of control law between the smooth transition of system state and the energy cost.

\begin{figure}[htbp]
\centering
\includegraphics[scale=0.30]{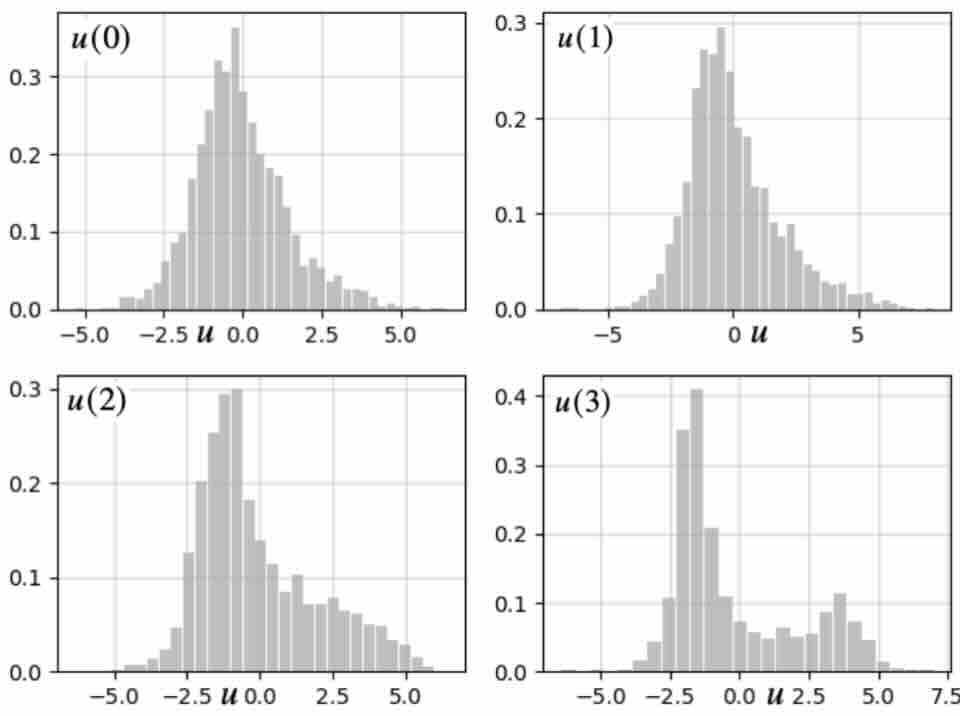}
\centering
\caption{The histograms of $u_{i}(k)$ at time step $k$ for each agent $i$ by cost function \eqref{costfunc1}. The upper left and right figures are $u_{i}(0)$ and $u_{i}(1), i = 1, \cdots, 1000$ respectively. The lower left and right figures are $u_{i}(2)$ and $u_{i}(3)$ respectively.}
\label{figd1}
\end{figure}

\begin{figure}[htbp]
\centering
\includegraphics[scale=0.30]{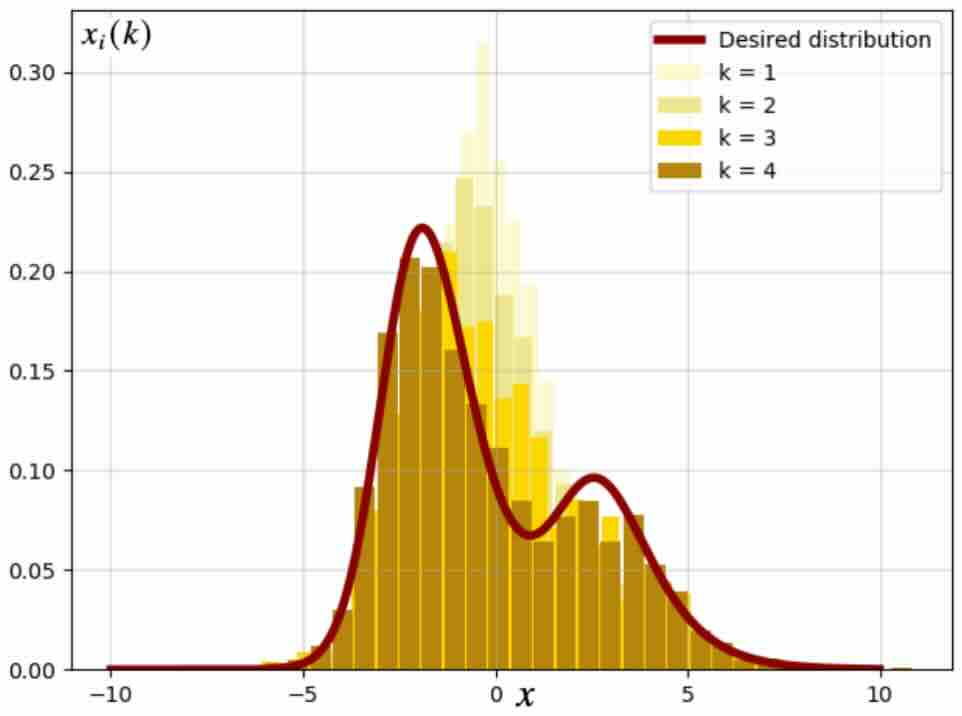}
\centering
\caption{The histograms of the system states $x_{i}(k)$ for $i = 1, \cdots, 2000$ at time steps $k = 1, 2, 3, 4$ by cost function \eqref{costfunc1}. The histogram at $K = 4$ is close to the specified terminal distribution \eqref{qt1}.}
\label{figd2}
\end{figure}

\begin{figure}[htbp]
\centering
\includegraphics[scale=0.30]{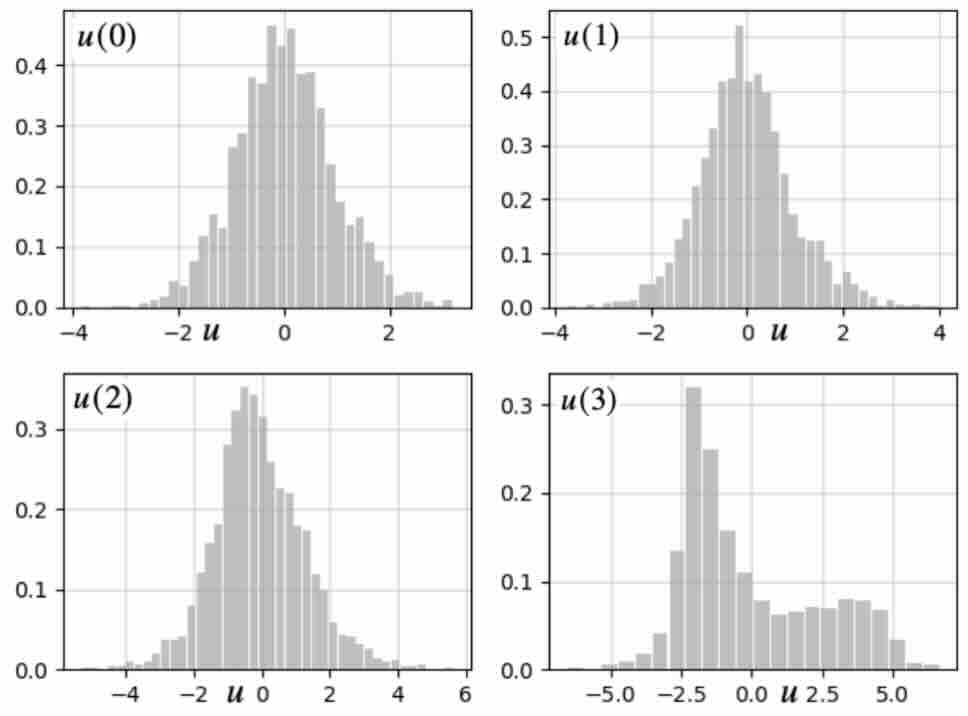}
\centering
\caption{The histograms of $u_{i}(k)$ at time step $k$ for each agent $i$ by cost function \eqref{WeightedCost}.}
\label{figd3}
\end{figure}

\begin{figure}[htbp]
\centering
\includegraphics[scale=0.30]{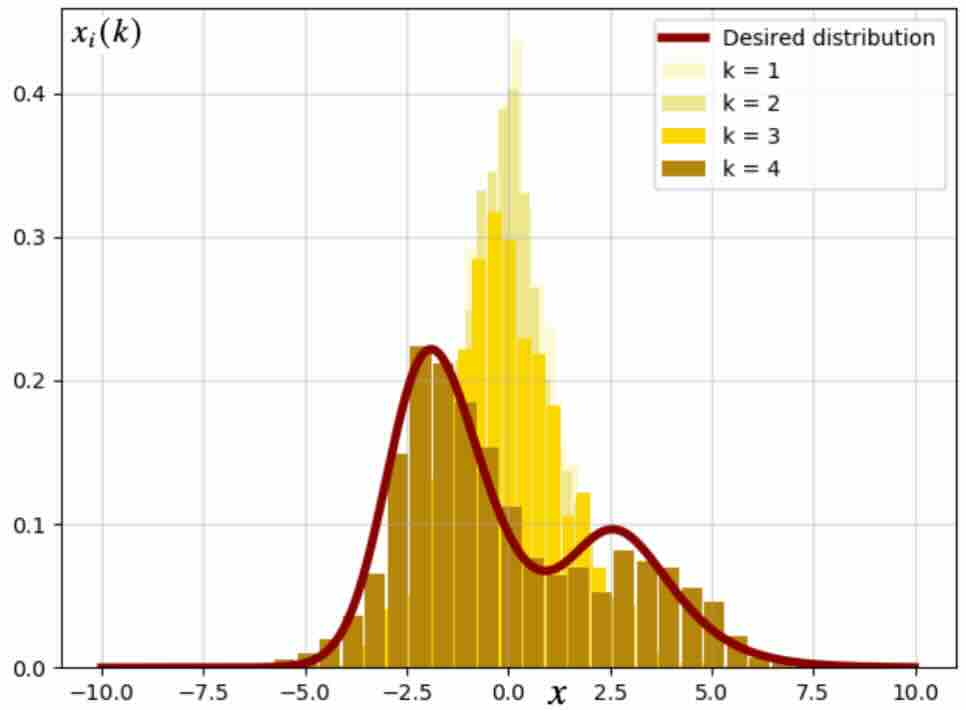}
\centering
\caption{The histograms of the system states $x_{i}(k)$ for $i = 1, \cdots, 2000$ at time steps $k = 1, 2, 3, 4$ by cost function \eqref{WeightedCost}. The histogram at $K = 4$ is close to the specified terminal distribution \eqref{qt1}.}
\label{figd4}
\end{figure}

Then we treat the discrete distribution (occupation measure) steering problem. The initial occupation measure $d\chi_{0}(\x)$ composes of the i.i.d. samples drawn from the the continuous distribution $d\chi_{0}(\x)$. Figure \ref{figd1} shows the histograms of the $u_{i}(k)$ for each agent at time step $k = 0, \cdots, 3$, by cost function \eqref{costfunc1}. Figure \ref{figd2} shows the histogram of the terminal occupation measure of the agents. The two peaks of the desired terminal state, of which the distribution is a mixture of two generalized logistic distributions, are well located at the desired points $x = -3$ and $x = 2$. The histogram in Figure \ref{figd2} is very close to $\chi_{T}(\x)$ in \eqref{qt1}, which validates the performance of our proposed algorithm. 

For the cost function of weighted energy effort and system energy \eqref{WeightedCost}, the histograms of the control inputs $u_{i}(k)$ are given in Figure \ref{figd3}. And the histogram of the terminal state of each agent $x_{i}(K)$ for $K = 4$ is shown in Figure \ref{figd4}, which is very close to the desired terminal distribution \eqref{qt1}. The results of two discrete distribution steering examples validate our proposed algorithm.

\subsection{Two Laplacians to two Gaussians}

Next, we consider the problem of steering the agents which are in separate groups to desired terminal groups. In this section, we simulate on steering a mixture of two Laplacians to a mixture of two Gaussians. Both initial and terminal distributions have two modes. The initial one is chosen as
\begin{equation}
\chi_0(\x)=\frac{0.5}{2} e^{|x-3|}+\frac{0.5}{2} e^{-|x+1|}
\label{q02}
\end{equation}
and the terminal one is specified as
\begin{equation}
\chi_T(\x)= \frac{0.5}{\sqrt{2 \pi}} e^{\frac{(x - 3)^2}{2}} + \frac{0.5}{\sqrt{2 \pi}} e^{\frac{(x + 3)^2}{2}}.
\label{qt2}
\end{equation}

The system parameters $a(k), k=0, \cdots, 3$ are i.i.d. samples drawn from the uniform distribution $U[0.3,0.5]$. The dimension of each $\mathfrak{U}(k)$ is $4$. 

\begin{figure}[htbp]
\centering
\includegraphics[scale=0.30]{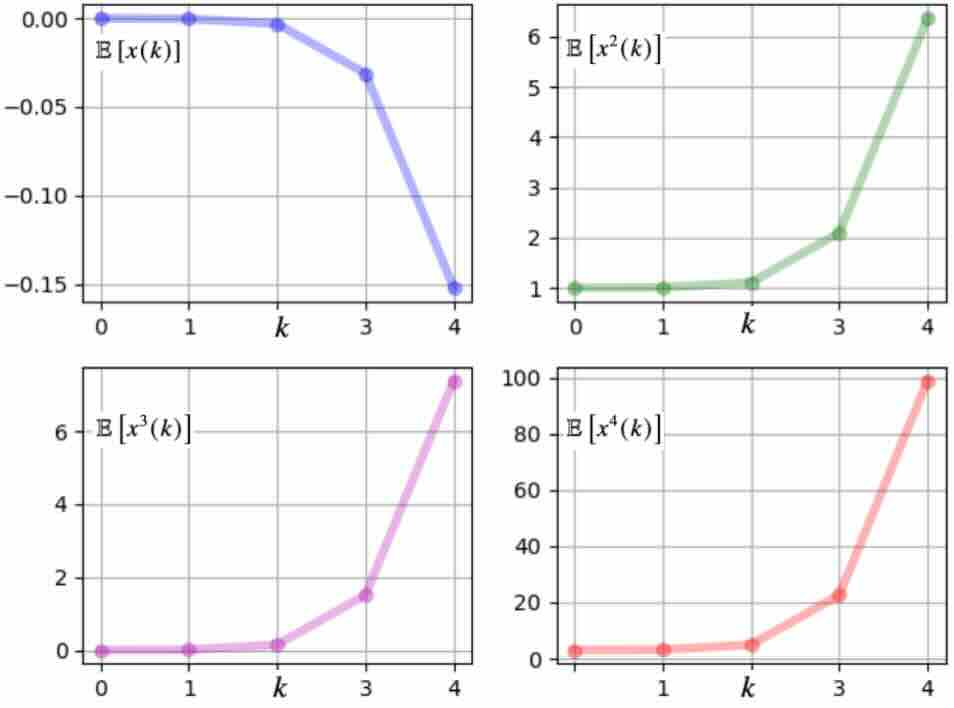}
\centering
\caption{$\mathfrak{X}(k)$ at time steps $k = 0, 1, 2, 3, 4$ by cost function \eqref{costfunc1}. The upper left figure shows $\mathbb{E}\left[ x(k)\right]$. The upper right one shows $\mathbb{E}\left[ x^{2}(k)\right]$. The lower left one shows $\mathbb{E}\left[ x^{3}(k)\right]$ and the lower right one shows $\mathbb{E}\left[ x^{4}(k)\right]$.}
\label{fig10}
\end{figure}

\begin{figure}[htbp]
\centering
\includegraphics[scale=0.30]{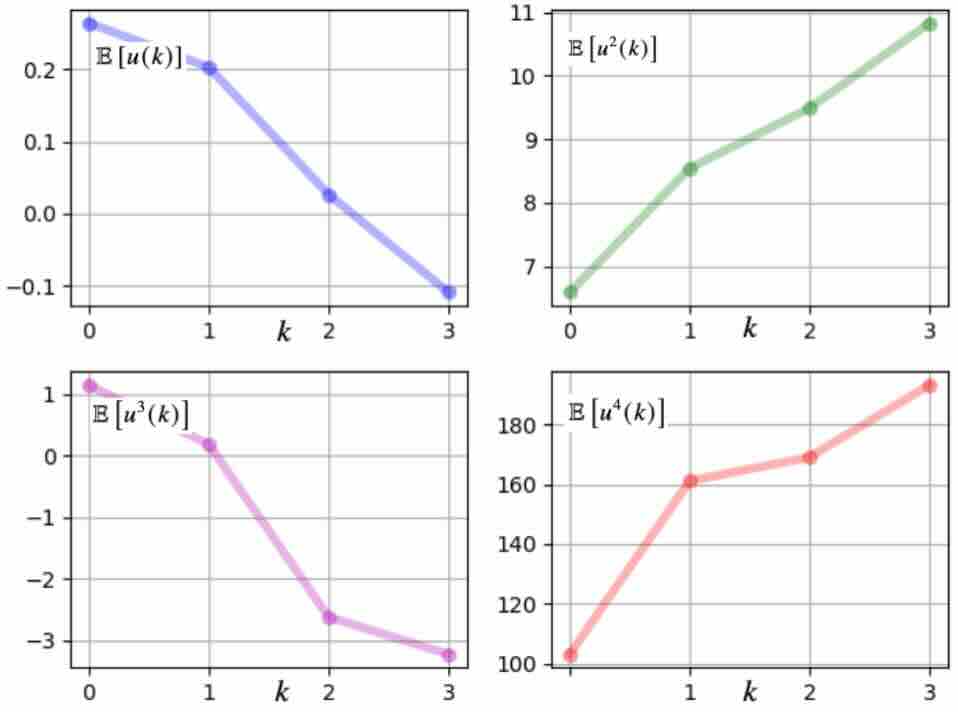}
\centering
\caption{$\mathfrak{U}(k)$ at time steps $k = 0, 1, 2, 3$ by cost function \eqref{costfunc1}. The upper left figure shows $\mathbb{E}\left[ u(k)\right]$. The upper right one shows $\mathbb{E}\left[ u^{2}(k)\right]$. The lower left one shows $\mathbb{E}\left[ u^{3}(k)\right]$ and the lower right one shows $\mathbb{E}\left[ u^{4}(k)\right]$.}
\label{fig11}
\end{figure}

\begin{figure}[htbp]
\centering
\includegraphics[scale=0.30]{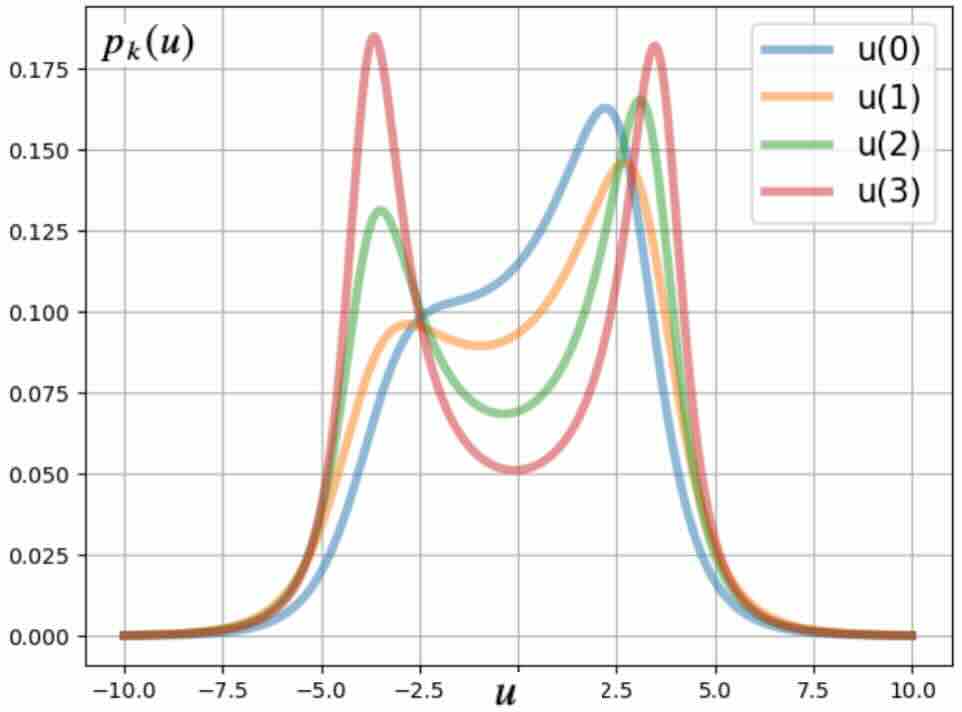}
\centering
\caption{Realized control inputs $\nu_{k}(\u)$ of $u(k)$ by $\mathfrak{U}(k)$ for $k = 0, 1, 2, 3$, which are obtained by cost function \eqref{costfunc1}.}
\label{fig12}
\end{figure}

We first perform the control task with the cost function \eqref{costfunc1}. The states of the moment system, i.e., $\mathscr{X}(k)$ for $k=0,1,2,3$, are illustrated in Figure \ref{fig10}. The controls of the moment system, i.e., $\mathfrak{U}(k)$ for $k=0,1,2,3$ are given in Figure \ref{fig11}. The realized controls in Figure \ref{fig12} also show that the transition of the control inputs is smooth, even the task is to steer a distribution with two modes to another one with two modes. The results of discrete distribution (occupation measure) steering is given in Figure \ref{figd5} and \ref{figd6}. The terminal distribution by the proposed steering scheme is close to the desired one. The energy effort $\sum_{k=0}^{3}\mathbb{E}\left[ u^{2}(k) \right] = 35.641$.

\begin{figure}[htbp]
\centering
\includegraphics[scale=0.30]{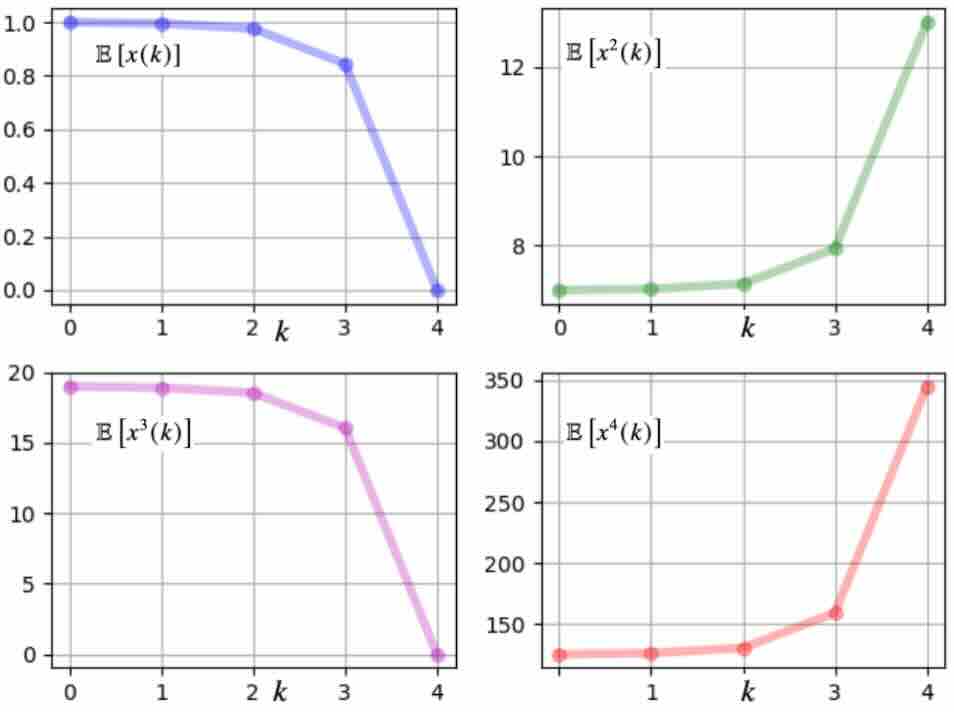}
\centering
\caption{$\mathfrak{X}(k)$ at time steps $k = 0, 1, 2, 3, 4$ by cost function \eqref{WeightedCost}.}
\label{fig13}
\end{figure}

\begin{figure}[htbp]
\centering
\includegraphics[scale=0.30]{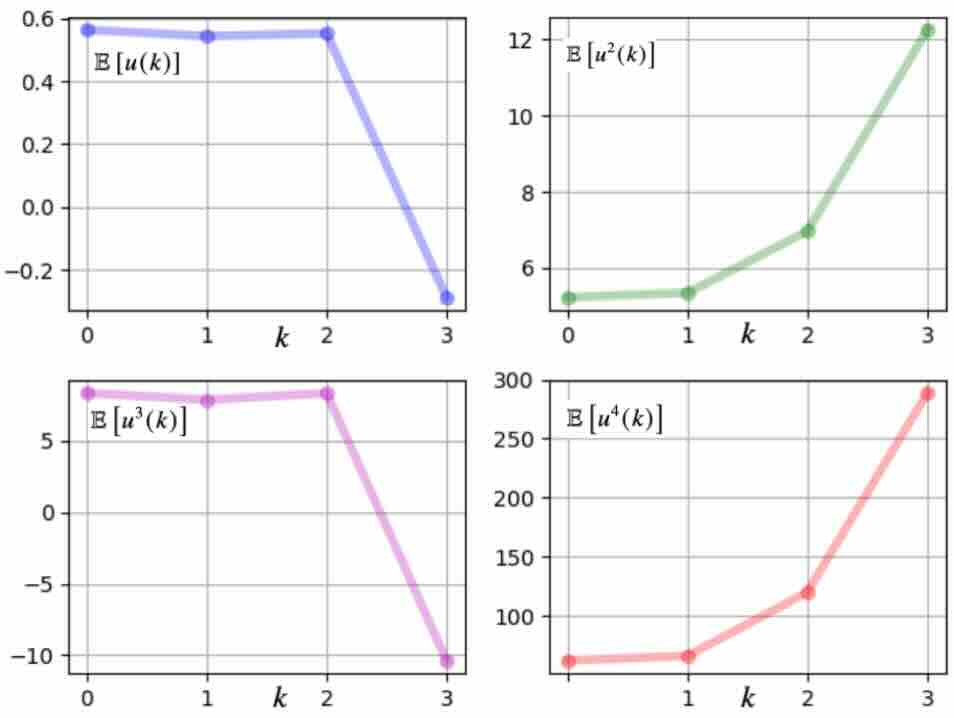}
\centering
\caption{$\mathfrak{U}(k)$ at time steps $k = 0, 1, 2, 3$ by cost function \eqref{WeightedCost}.}
\label{fig14}
\end{figure}

\begin{figure}[htbp]
\centering
\includegraphics[scale=0.30]{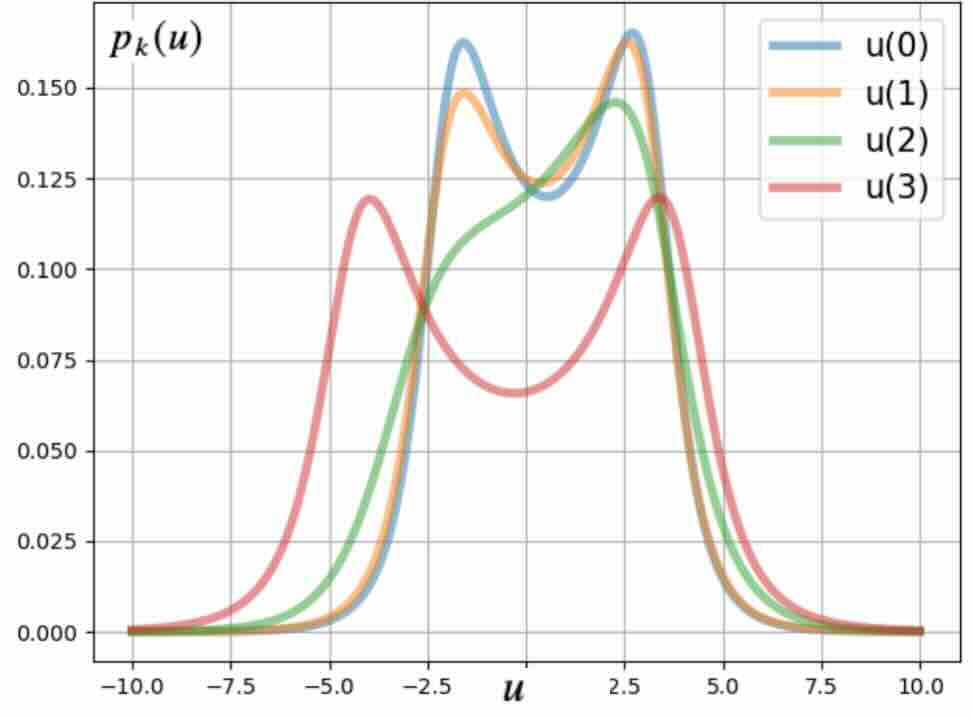}
\centering
\caption{Realized control inputs $\nu_{k}(u)$ of $u(k)$ by $\mathfrak{U}(k)$ for $k = 0, 1, 2, 3$, which are obtained by cost function \eqref{WeightedCost}.}
\label{fig15}
\end{figure}

\begin{figure}[htbp]
\centering
\includegraphics[scale=0.30]{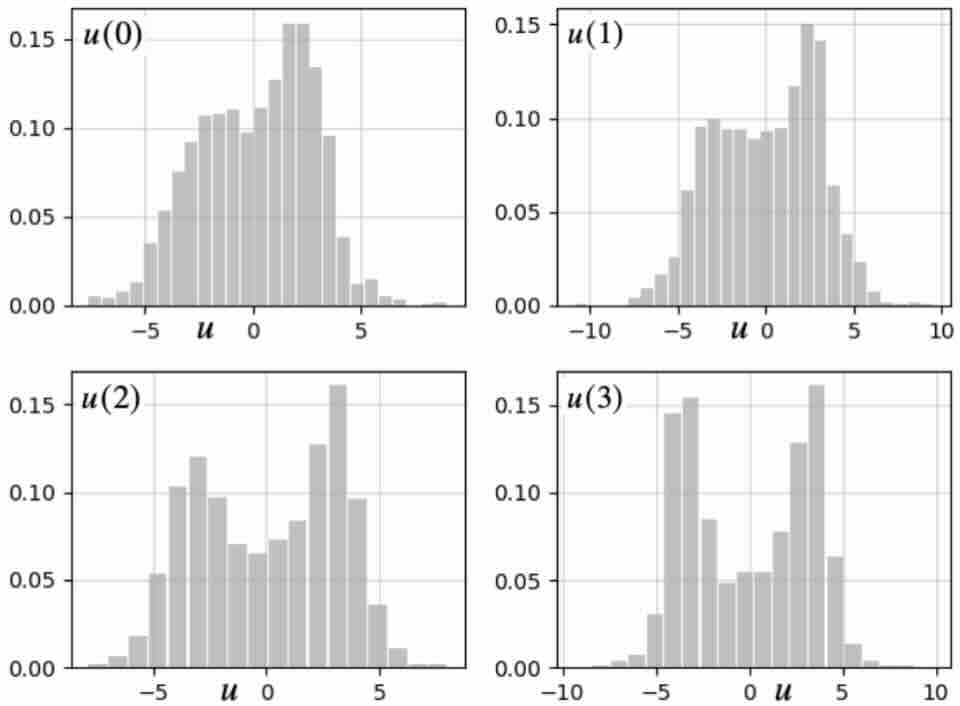}
\centering
\caption{The histograms of $u_{i}(k)$ at time step $k$ for each agent $i$ by cost function \eqref{costfunc1}. The upper left and right figures are $u_{i}(0)$ and $u_{i}(1), i = 1, \cdots, 1000$ respectively. The lower left and right figures are $u_{i}(2)$ and $u_{i}(3)$ respectively.}
\label{figd5}
\end{figure}

\begin{figure}[htbp]
\centering
\includegraphics[scale=0.30]{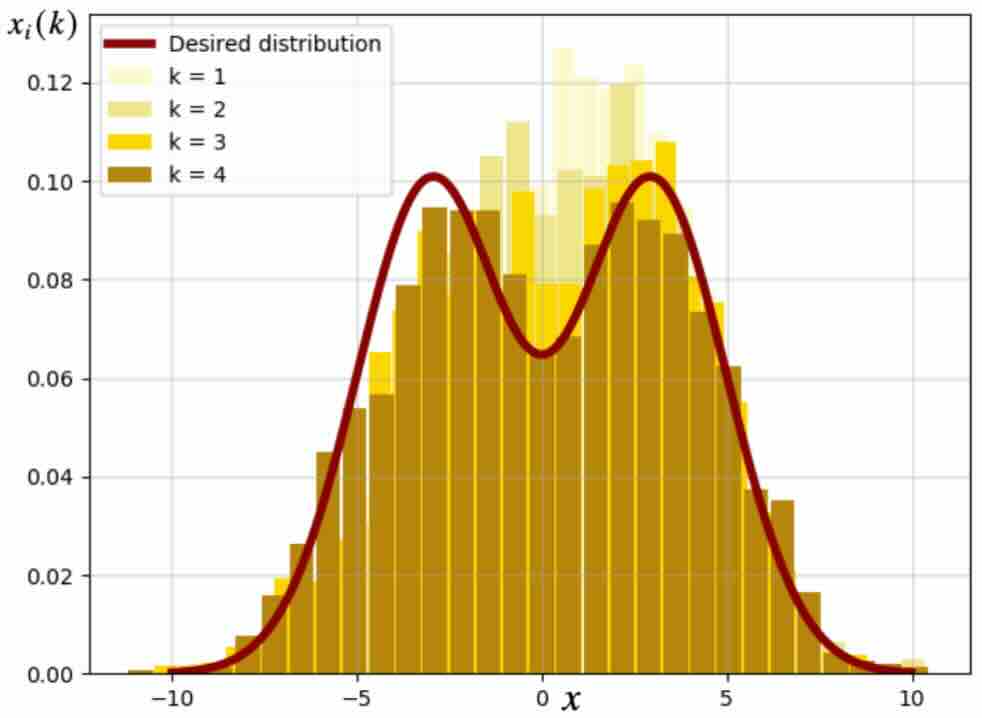}
\centering
\caption{The histograms of the system states $x_{i}(k)$ for $i = 1, \cdots, 2000$ at time steps $k = 1, 2, 3, 4$ by cost function \eqref{costfunc1}. The histogram at $K = 4$ is close to the specified terminal distribution \eqref{qt1}.}
\label{figd6}
\end{figure}

\begin{figure}[htbp]
\centering
\includegraphics[scale=0.30]{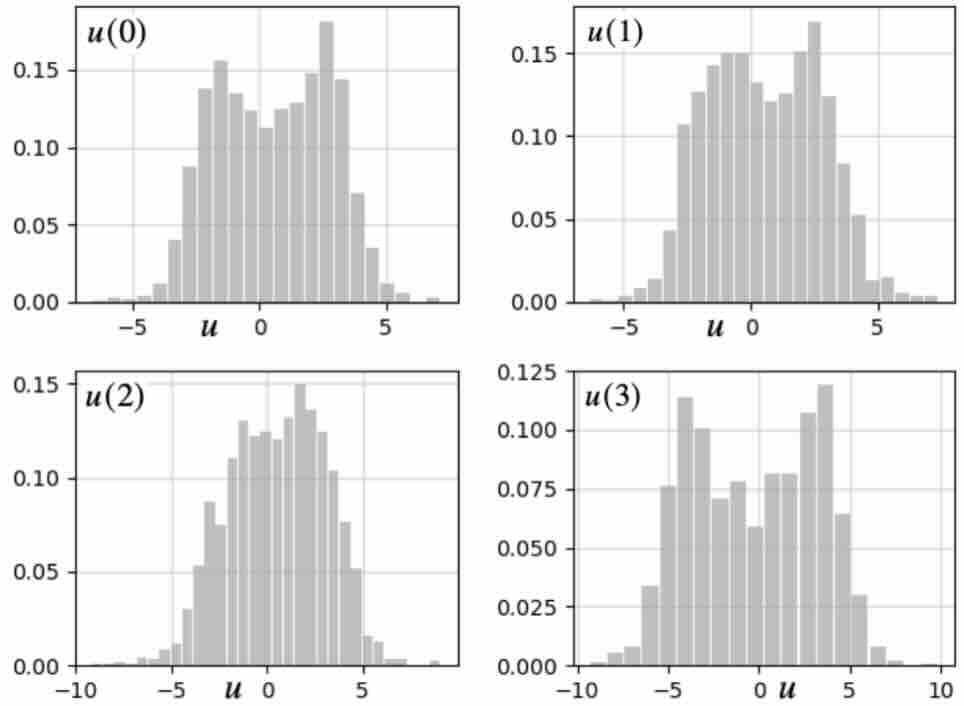}
\centering
\caption{The histograms of $u_{i}(k)$ at time step $k$ for each agent $i$ by cost function \eqref{WeightedCost}.}
\label{figd7}
\end{figure}

\begin{figure}[t]
\centering
\includegraphics[scale=0.30]{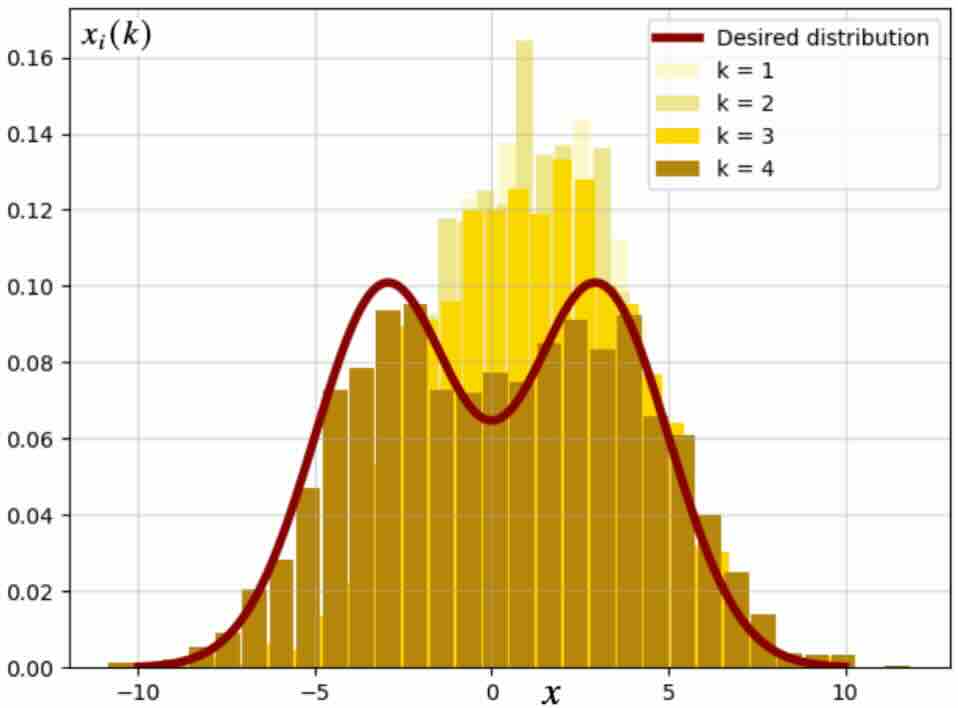}
\centering
\caption{The histograms of the system states $x_{i}(k)$ for $i = 1, \cdots, 2000$ at time steps $k = 1, 2, 3, 4$ by cost function \eqref{WeightedCost}. The histogram at $K = 4$ is close to the specified terminal distribution \eqref{qt2}.}
\label{figd8}
\end{figure}

Next, we do optimization \eqref{optimization1} with the cost function \eqref{costfunc3}. In this simulation, we choose the cost function as \eqref{WeightedCost}. The simulation results are given in Figure \ref{fig13}, \ref{fig14} and \ref{fig15}. The results of discrete distribution (occupation measure) steering is given in Figure \ref{figd7} and \ref{figd8}. The histogram of the terminal states of the agents is close to the desired continuous terminal distribution \eqref{qt2}, which reveals the performance of our proposed algorithm. The energy effort $\sum_{k=0}^{3}\mathbb{E}\left[ u^{2}(k) \right] = 29.783$. 

\section{A concluding remark}
We consider the general distribution steering problem where the distributions to steer are arbitrary, which are only required to have first several orders of finite power moments. In our previous paper \cite{wu2023density}, we proposed a moment counterpart of the primal system for control. However, we were not able to put forward a control law based on optimization in the manner of conventional optimal control, which makes it hardly possible for us to obtain the optimal control inputs by specific purposes, such as minimum energy effort. In this paper, we investigate the general distribution steering problem by optimization. The domain of the control inputs of the moment system is not convex and has a complex topology, which causes difficulty in optimization. We prove the controllability of the moment system and propose a subset as the domain for optimization of which the convexity is proved. The subset is proved to be a set of solutions minimizing a weighted sum of integrated squared distances. Then we consider different types of cost functions, including the smoothness of the state transition, the system energy, the control energy effort, and a general form of cost function. A realization of the control inputs by the squared Hellinger distance is given to put forward a control scheme for the general distribution steering problem. We consider two problems in swarm robotics and formulate them as two distribution steering problems for simulation. The numerical results of the simulations validate our proposed algorithms. By the simulation results, we note that to yield smooth transition of the system states, one may need more control energy effort. Moreover, moment approximation provides a much more computationally efficient approach to distribution steering compared to discretizing the support of the distribution, revealed by the results of this paper.

In the future work, we would like to extend the results of this paper to nonlinear systems. We would also like to extend the results of the first-order system to more general systems, which will not be a trivial extension since the positive definiteness of the Hankel matrix will no longer be the sufficient and necessary condition for the existence of the multi-dimensional control inputs. Extending numerous results of this paper to the multidimensional systems is a difficult task which will require mathematical tools from real algebraic geometry and other subjects. 

\bibliographystyle{plain}
\bibliography{autosam}

\begin{wrapfigure}{l}{20mm} 
\includegraphics[width=1in,height=1.25in,clip,keepaspectratio]{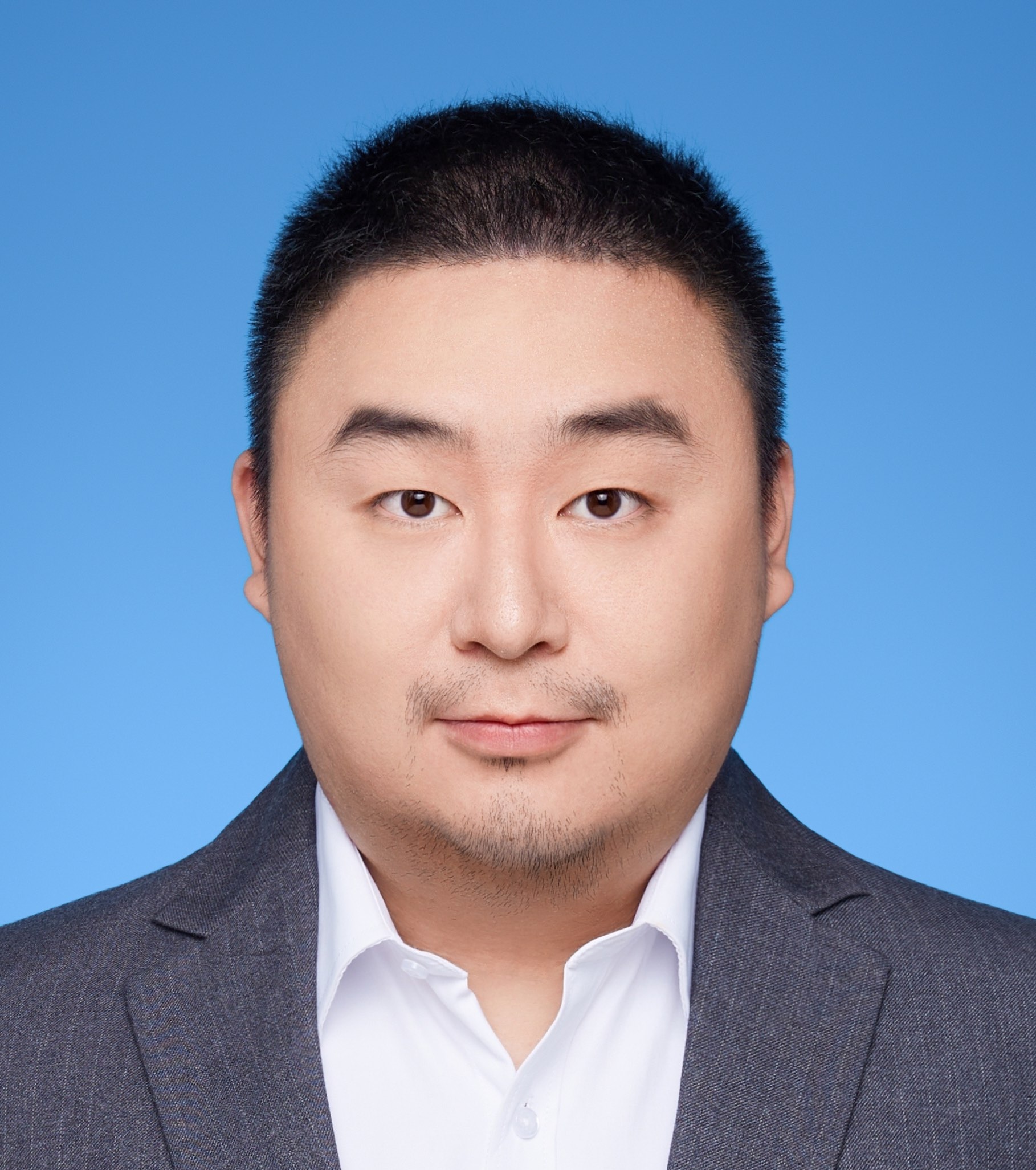}
  \end{wrapfigure}\par
  \textbf{Guangyu Wu} received the Ph.D. degree in Control Science and Engineering (advisor: Anders Lindquist) from Shanghai Jiao Tong University, Shanghai, China, in 2024.

He is currently a research scientist with the College of Computing and Data Science, Nanyang Technological University, Singapore. His current research interests encompass the broad fundamentals of digital twins for the data center, including system modelling and reduction, machine learning (with a focus on fine-resolution and few-shot surrogate modelling in AI4Science), prediction, and the optimization and control. He is a recipient of the Eric and Wendy Schmidt AI in Science Postdoctoral Fellowship. He serves as an active reviewer for IEEE Transactions on Automatic Control, Automatica, IEEE Control Systems Letters, and the IEEE Conference on Decision and Control. 

\begin{wrapfigure}{l}{25mm} 
\includegraphics[width=1in,height=1.25in,clip,keepaspectratio]{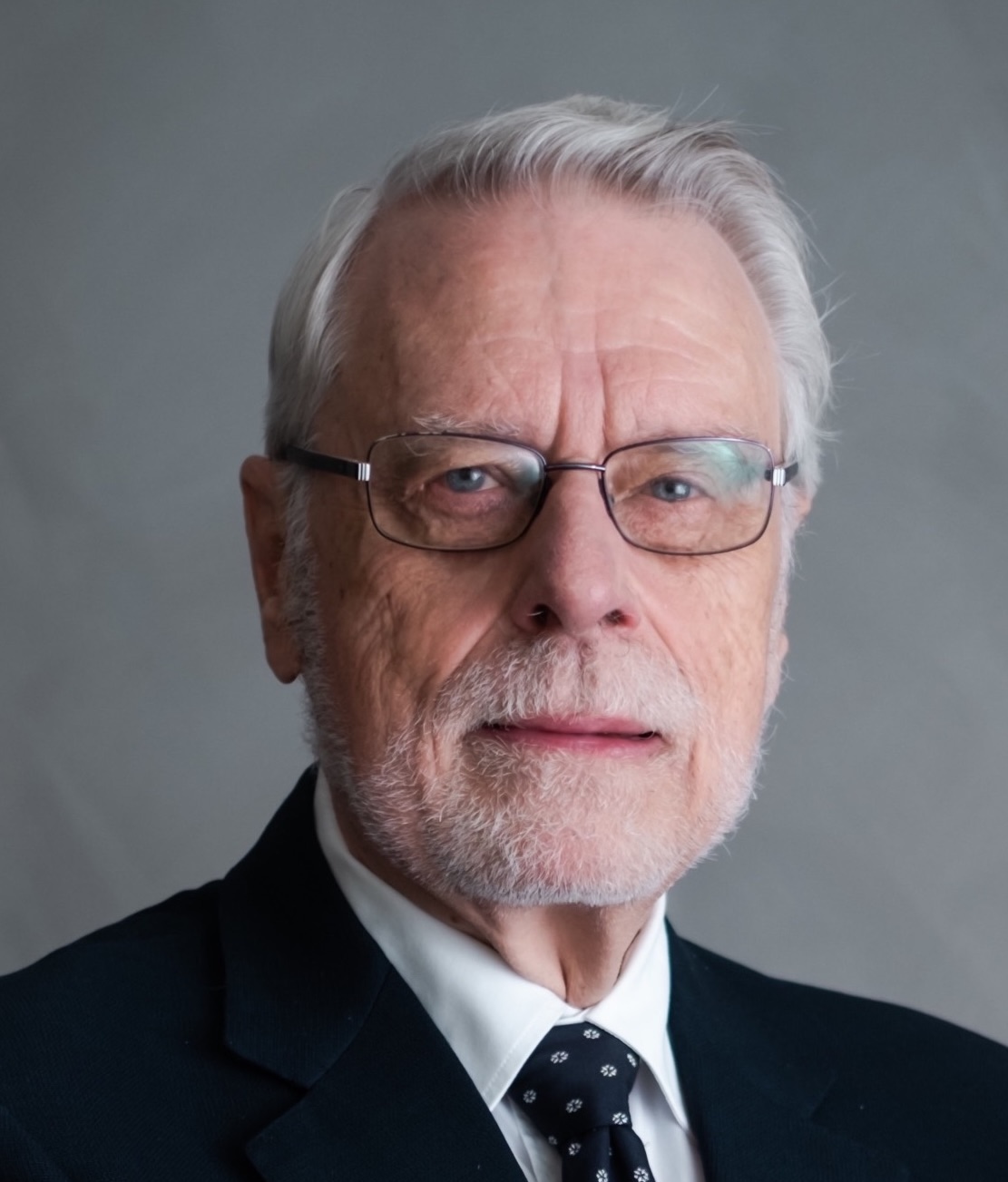}
  \end{wrapfigure}\par
  \textbf{Anders Lindquist} received the Ph.D. degree in Optimization and Systems Theory from the Royal Institute of Technology (KTH), Stockholm, Sweden, in 1972, an honorary doctorate (Doctor Scientiarum Honoris Causa) from Technion (Israel Institute of Technology) in 2010 and Doctor Jubilaris from KTH in 2022.

He is currently a Distinguished Professor at Anhui University, Hefei, China, Professor Emeritus at Shanghai Jiao Tong University, China, and Professor Emeritus at the Royal Institute of Technology (KTH), Stockholm, Sweden. Before that he had a full academic career in the United States, after which he was appointed to the Chair of Optimization and Systems at KTH.

Dr. Lindquist is a Member of the Royal Swedish Academy of Engineering Sciences, a Foreign Member of the Chinese Academy of Sciences, a Foreign Member of the Russian Academy of Natural Sciences (elected 1997), a Member of Academia Europaea (Academy of Europe), an Honorary Member the Hungarian Operations Research Society, a Life Fellow of IEEE, a Fellow of SIAM, and a Fellow of IFAC. He received the 2003 George S. Axelby Outstanding Paper Award, the 2009 Reid Prize in Mathematics from SIAM, and the 2020 IEEE Control Systems Award, the IEEE field award in Systems and Control.

\end{document}